\documentclass{amsart}
\usepackage{amssymb,latexsym,enumerate,amscd}
\newtheorem{Thm}{Theorem}[section]
\newtheorem{Prop}[Thm]{Proposition}
\newtheorem{Lem}[Thm]{Lemma}
\newtheorem{Cor}[Thm]{Corollary}

\newtheorem{Ques}[Thm]{Question}
\newtheorem{Add}[Thm]{Addendum}
\theoremstyle{remark}
\newtheorem{Rem}[Thm]{Remark}
\theoremstyle{definition}
\newtheorem*{Not}{Notation}
\newtheorem*{Def}{Definition}
\numberwithin{equation}{section}
\DeclareMathOperator{\im}{im}
\DeclareMathOperator{\id}{id}
\DeclareMathOperator{\GL}{GL}
\DeclareMathOperator{\Wh}{Wh}

\DeclareMathOperator{\Iso}{Iso}
\begin{document}

\title[$K$-theory of solvable groups]
{$K$-theory of solvable groups}

\author[F. T. Farrell]{F. Thomas Farrell}
\address{Department of Mathematical Sciences \\
State University of New York at Binghamton \\
Binghamton \\
NY 13902-6000 \\
USA}
\email{farrell@math.binghamton.edu}
\urladdr{http://www.math.binghamtom.edu/farrell/}
\thanks{The first author was supported in part by the National
Science Foundation}

\author[P. A. Linnell]{Peter A. Linnell}
\address{Department of Mathematics \\
Virginia Tech \\
Blacksburg \\
VA 24061-0123 \\
USA}
\email{linnell@math.vt.edu}
\urladdr{http://www.math.vt.edu/people/linnell/}

\begin{abstract}
We first prove that the Whitehead group of a torsion-free
virtually solvable linear group vanishes.  Next we
make a reduction of the fibered isomorphism conjecture from 
virtually solvable groups to a class of virtually solvable $\mathbb
{Q}$-linear groups.  Finally we prove an $L$-theory analogue for
elementary amenable groups.
\end{abstract}

\keywords{solvable group, Whitehead group, fibered isomorphism
conjecture}

\subjclass{Primary: 19B28; Secondary: 19G24}

\date{%
Mon Jul 15 18:23:40 EDT 2002}

\maketitle

\section{Introduction} \label{Sintroduction}

This paper will be concerned about verifying some
isomorphism conjectures in $K$-theory for solvable groups.  Recall
that if $\mathcal {C}$ is a class of groups, then a group is
virtually $\mathcal {C}$ if it has a $\mathcal {C}$-subgroup of
finite index (equivalently it has a normal
$\mathcal {C}$-subgroup of
finite index).  Let $\Wh(G)$ denote the Whitehead group of the group
$G$; that is $K_1(\mathbb {Z}G)/\{\pm G\}$ and if $R$ is any ring,
let $\GL_n(R)$ denote the group of $n$ by $n$ invertible matrices
over $R$.  Our first result is
\begin{Thm} \label {T1}
Let $G$ be a torsion-free virtually solvable
subgroup of $\GL_n(\mathbb {C})$.  Then $\Wh (G) = 0$.
\end{Thm}
Let us say that $G$ is a linear group means that $G$ is isomorphic
to a subgroup of $\GL_n(\mathbb {C})$.  Then
another way of stating Theorem \ref{T1} is that if $G$ is a
torsion-free virtually solvable linear group, then $\Wh (G) =
0$.  Of course if $G$ is a torsion-free virtually solvable linear
group, then so is $G \times \mathbb {Z}^n$ for all integers $n \ge
0$.  Let $\tilde{K}_0(\mathbb {Z}G)$ denote the reduced projective
class group, i.e.\ $K_0(\mathbb {Z}G)
/\langle [\mathbb {Z}G] \rangle$.
It now follows from Bass's contracted functor argument
\cite[chapter XII, \S 7]{bass} that $\tilde {K}_0(\mathbb {Z}G)
= 0$ and $K_i(\mathbb {Z}G) = 0$ for all negative
integers $i$.

Our second result concerns the Fibered Isomorphism Conjecture (FIC)
\cite{farrelljonesiso}.  A precise formulation of this conjecture is
given at the beginning of appendix, Section \ref{Sappendix}.
It is formulated there for any homotopy invariant
spectrum valued functor $\mathcal {S}$ defined on the
category of topological spaces.  However it is only plausibly true
for certain functors; e.g.\ $\mathcal {P}$, $\mathcal {K}$ and
$\mathcal {L}^{-\infty}$.  When FIC is mentioned in Sections
\ref{Sintroduction}--\ref{ST2}, it always refers to this formulation
where $\mathcal {S}$ is the stable topological psuedo-isotopy
functor $\mathcal {P}$.
For the purposes of this paper let us make
the following definition (in this paper
$\otimes$ will always mean $\otimes_{\mathbb{Z}}$).
\begin{Def}
Let $G$ be a group.  Then $G$ is \emph{nearly crystallographic}
means that $G$ is finitely generated and
there exists $A \lhd G$ such that $A$ is torsion-free abelian
of finite rank (i.e.\ is isomorphic to a subgroup of $\mathbb
{Q}^n$ for some integer $n$), $C \leqslant G$ such that $C$ is
virtually cyclic, $A \cap C = 1$ and $AC = G$ (so $G = A
\rtimes C$), and the conjugation action of $C$ on $A$ makes
$A\otimes \mathbb {Q}$ into an
irreducible $\mathbb {Q}C$-module.
\end{Def}
It is easy to see that nearly crystallographic groups are rationally
linear, i.e.\ isomorphic to subgroups of $\GL_n(\mathbb {Q})$.  From
\cite[theorem 2.1 and theorem A.8]{farrelljonesiso}, we see that
the FIC is true for subgroups of cocompact discrete subgroups of
virtually connected Lie groups.
Also \cite[theorem 4.8]{farrelljonesiso},
shows that the FIC is true for nearly
crystallographic groups in the case $A = \mathbb {Z}^n$.
On the other hand the FIC remains unproven for the case $G = \mathbb
{Z}[\frac{1}{2}] \rtimes_{\alpha} \mathbb {Z}$, where $\alpha$ is
multiplication by 2.  We shall prove
\begin{Thm} \label{T2}
Suppose the FIC is true for all nearly crystallographic
groups.  Then the FIC is true for all virtually solvable groups.
\end{Thm}
In the course of the proof of Theorem \ref{T2}, we establish the FIC
for the group $(A \wr \mathbb {Z}) \wr F$ where $F$ is a finite
group and $A$ is either an elementary abelian $p$-group
for some prime $p$ or a torsion-free abelian
group (without the assumption that the FIC is true for nearly
crystallographic groups).
This statement follows immediately from Lemma
\ref{Lwreathp}, Lemma \ref{Lwreath} and Theorem \ref{Tdirlimit}.

In Section \ref{SLtheory} we shall prove an analogue of the FIC for
$L$-theory of elementary amenable groups and we shall apply this
result in Section \ref{SLtheoryapplication}.
Finally in the appendix, Section \ref{Sappendix} at
the end, we verify that the FIC remains true under direct limits.
This result is crucial for Sections \ref{ST2} and \ref{SLtheory}.

This work was carried out while we were at the
Sonderforschungsbereich in M\"unster.  We would like to thank
Wolfgang L\"uck for organizing our visits to M\"unster, and the
Sonderforschungsbereich for financial support.
We also wish to thank the referee for
pointing out an error in the preprint of this article
(see Preprintreihe SFB 478, Heft 103 at the Universit\"at
M\"unster) where we asserted that the bottom horizontal arrow in
diagram \eqref{EAcomlim} is a monomorphism.

\section{The FIC and related required results}

In this section we will review some mostly known results related to
the FIC.

\begin{Not}
If $p$ is a positive integer, then $\mathbb {Z}_p$ will indicate
the cyclic group of order $p$.  The direct sum of $n$ copies of the
group $A$ will be indicated by $A^n$.  If $G$ is a group acting
faithfully on a set $X$, then $A \wr_X G$ will denote the wreath
product with respect to this action
(see \cite[\S 2.6]{dixonmortimer} or \cite[\S 1.6]{robinson}).
With regard to the situation $X$ infinite, this will always mean
the restricted wreath product, so the base group is $\bigoplus_{x
\in X} A$ (as opposed to the Cartesian product).
In the case $X = G$ and $G$ is acting
regularly on $X$, we shall write $A \wr G$ in
place of $A \wr_X G$.  
\end{Not}

\begin{Lem} \label{Lficwh}
Let $H\lhd G$ be groups such that $G/H$ is finitely generated
virtually torsion-free abelian, and let $p \colon G
\twoheadrightarrow G/H$ denote
the natural epimorphism.  If $\Wh (p^{-1}(S) \times \mathbb {Z}^n) =
0$ for all virtually cyclic subgroups $S$ of $G/H$ and for all
non-negative integers $n$, then $\Wh(G) = 0$.
\end{Lem}
\begin{proof}
This follows from \cite[proposition 2.4]{farrelljonesiso}.
\end{proof}

\begin{Lem} \label{Lficinduction}
Let $p \colon \Gamma \twoheadrightarrow G$ be an epimorphism of
groups such that the FIC is true for $G$ and
also for $p^{-1}(S)$ for all
virtually cyclic subgroups $S$ of $G$.  Then the FIC is also true
for $\Gamma$.
\end{Lem}
\begin{proof}
\cite[proposition 2.2]{farrelljonesiso}.
\end{proof}

\begin{Lem} \label{Lficsubgroup}
Let $\Gamma$ be a group.  If the FIC is true for $\Gamma$,
then the FIC is true for every subgroup of $\Gamma$.
\end{Lem}
\begin{proof}
\cite[theorem A.8]{farrelljonesiso}
\end{proof}

\begin{Lem} \label{Lficfree}
Let $F$ be a finitely generated free group, let $X$ be a finite
set, and let $S$ denote the symmetric group on $X$.
Then the FIC is true for $F\wr_X S$.
\end{Lem}
\begin{proof}
This is really \cite[fact 3.1]{farrellroushon}; we sketch the
argument.  Since $F$ is finitely generated free, it
is isomorphic to a subgroup of a
cocompact lattice $G$ in $\Iso (\mathbb {H}^2)$ (the group of
isometries of hyperbolic 2-space) and
then $F\wr_X S$ is isomorphic to a
subgroup of $G \wr_X S$.  By Lemma \ref{Lficsubgroup},
it will be sufficient to show that $G \wr_X S$ satisfies the FIC.
Let $n = |X|$.  Then $G^n$ and hence also $G \wr_X S$ are
cocompact lattices in $\Iso (\mathbb {H}^2
\times \dots \times \mathbb {H}^2)$ (where there are $n$ copies of
$\mathbb {H}^2$ in the direct product).  The result now follows from
\cite[proposition 2.3]{farrelljonesiso}.
\end{proof}

Two useful results in conjunction with Lemma \ref{Lficfree}
are the following.

\begin{Lem} \label{Lwreathwreath}
Let $A,B,G$ be groups, and let $X,Y$ be finite sets.  Suppose $A$
acts on $X$ faithfully and $B$ acts on $Y$ faithfully.  Then $(G
\wr_X A) \wr_Y B \cong G \wr_Z C$ for some group $C$ and some finite
set $Z$ with $C$ acting faithfully.
\end{Lem}
\begin{proof}
The isomorphism of \cite[exercise 2.6.7]{dixonmortimer}
or \cite[1.6.4(ii)]{robinson} tells us that
\[
(G \wr_X A) \wr_Y B \cong G \wr_{X \times Y} (A\wr_Y B)
\]
for some suitable action of $A \wr_Y B$ on $X \times Y$.
It is easily checked
that this action is faithful, so the result follows with $Z = X
\times Y$ and $C = A \wr_Y B$.
\end{proof}

\begin{Lem} \label{Lalgebraiclemma}
Let $n$ be a positive integer, let $F$ be a group of order $n$,
and let $H \lhd G$ be groups with $|G/H| = n$.  Then $H \wr F$ has a
subgroup isomorphic to $H^n$ and $H \wr G/H$ has a subgroup
isomorphic to $G$.
\end{Lem}
\begin{proof}
The first statement is obvious, while the second follows from either
\cite[theorem 2.6A]{dixonmortimer} or
\cite[an algebraic lemma]{farrellroushon}.
\end{proof}

Finally we need

\begin{Lem} \label{Lficvirabelian}
Let $\Gamma$ be a virtually abelian group.  Then $\Gamma$ satisfies
the FIC.
\end{Lem}
\begin{proof}
By Theorem \ref{Tdirlimit}, we may assume that $\Gamma$ is finitely
generated.  The result now follows from \cite[proposition
2.4]{farrelljonesiso}.
\end{proof}

\section{Proof of Theorem \ref{T1}} \label{ST1}

\begin{proof}[Proof of Theorem \ref{T1}]
Since Whitehead groups commute with direct limits,
we may assume that $G$ is finitely generated (so in particular all
further groups considered in this proof can be assumed to be
countable).  By a theorem of Malcev
\cite[15.1.4]{robinson} we know that there exists $H \lhd A \lhd G$
with $H\lhd G$
such that $H$ is nilpotent, $A/H$ is abelian,
and $G/A$ is finite.  Since $G$ is finitely generated, we may assume
that $A/H$ is free abelian of finite rank.
Let $C/H$ be a virtually cyclic subgroup of
$G/H$, let $n$ be a non-negative integer, and set $D = C \times
\mathbb {Z}^n$.  By Lemma \ref{Lficwh}, it will
suffice to show that $\Wh (D) = 0$.  Set $E = H \times \mathbb
{Z}^n$.  Then $E$ is a normal nilpotent subgroup of $D$ and $D/E$ is
virtually cyclic.  Since $D/E$ is virtually cyclic, it has a finite
normal subgroup $F/E$ such that $D/F$ is either infinite cyclic or
infinite dihedral.  We now apply Waldhausen's results
\cite{waldhausen} and we shall use some of his terminology.
We shall also require the following well known result; see for
example
\cite[lemmas 23b and 24b]{farrellhsiangformula}
or \cite[p.~247]{waldhausen} and \cite[theorem
VIII.3.1]{brown}.
\begin{Lem} \label{Lregnilpotent}
Let $G$ be a finitely generated torsion-free virtually nilpotent
group.  Then $\mathbb {Z}G$ is regular Noetherian.
\end{Lem}
Since a countable group is the ascending union of
its finitely generated subgroups, we may apply \cite[theorem
19.1(iii)]{waldhausen} to deduce that $\mathbb {Z}F$ is regular
coherent.  We want to show that $\Wh (D) = 0$, and
we have two cases to consider.
\begin{enumerate}
\item
$D/F$ is infinite cyclic.  Since a finitely generated virtually
nilpotent group is virtually poly-$\mathbb {Z}$, we see from
\cite[theorem 3.2]{farrellhsiangjlms} that $\Wh (F_0) = 0$
for all finitely
generated subgroups $F_0$ of $F$,
consequently $\Wh (F) = 0$ because Whitehead groups
commute with direct with direct limits.  We can now apply
\cite[corollary 17.2.3]{waldhausen} to deduce that $\Wh (D) = 0$.

\item
$D/F$ is infinite dihedral.  Then there are subgroups $D_1/F, D_2/F$
of order two in $D/F$ such that  $D \cong D_1 *_F D_2$.  Also
$D_0$ is virtually poly-$\mathbb {Z}$ for all
finitely generated subgroups of $D_1$ and $D_2$.  Therefore
$\Wh(D_0) = 0$ by \cite[theorem 3.2]{farrellhsiangjlms}
and hence $\Wh (D_1)
=\Wh(D_2) = 0$, because Whitehead groups commute with direct limits.
We now deduce from \cite[corollary 17.1.3]{waldhausen} that $\Wh (D)
= 0$.
\end{enumerate}
It now follows from Lemma \ref{Lficwh} that $\Wh (G) =
0$, as required.
\end{proof}

\begin{Rem}
One could also consider linear groups over fields of nonzero
characteristic; i.e.\
subgroups of $\GL_n(k)$ for fields $k$ of nonzero characteristic.
However since
the groups we are considering in Theorem \ref{T1} are torsion-free,
\cite[15.1.1]{robinson} would tell us that such a group would have a
normal diagonalizable subgroup of finite index, consequently such a
group would be virtually abelian and so we would not be
getting anything new.
\end{Rem}

\section{Proof of Theorem \ref{T2}} \label{ST2}

\begin{Lem} \label{Lwreathp}
Let $n$ be a non-negative integer, let $p$ be a prime,
and let $F$ be a finite group.
Then $(\mathbb {Z}_p^n \wr \mathbb {Z}) \wr F$ satisfies the FIC.
\end{Lem}
\begin{proof}
Fix an abelian group $A$ and a positive integer $i$.
Let $G_i = A^{i+1}*_{A^i}$, the HNN extension with respect to the
following two embeddings of $A^i$ in $A^{i+1}$:
\[
(a_1, \dots, a_i) \mapsto (a_1, \dots, a_i,0)
\quad\text{and}\quad
(a_1, \dots, a_i) \mapsto (0, a_1, \dots, a_i).
\]
Then $A \wr \mathbb {Z} = \varinjlim G_i$ and hence
\[
(A \wr \mathbb {Z}) \wr F = \varinjlim (G_i \wr F).
\]
Consequently by Theorem \ref{Tdirlimit}, it suffices to show that
each $G_i \wr F$ satisfies the FIC when $A = \mathbb {Z}_p^n$.
Since $G_i$ is the fundamental group of a finite graph of groups
with finite vertex groups, we see from
\cite[theorem IV.1.6]{dicksdunwoody}
that $G_i$ has a normal finitely generated
free subgroup $H$ of finite index.
Using Lemma \ref{Lalgebraiclemma},
we may embed $G_i$ in $H \wr
G_i/H$, and then $G_i \wr F$ embeds in $(H \wr G_i/H) \wr
F$.  The result now follows from Lemmas \ref{Lficsubgroup},
\ref{Lficfree} and \ref{Lwreathwreath}.
\end{proof}

Using Lemma \ref{Lwreathp}, we obtain the following result.
\begin{Cor} \label{Cwreathp}
Let $\Gamma$ be a group which is torsion abelian by virtually
cyclic.  Then $\Gamma$ satisfies the FIC.
\end{Cor}
\begin{proof}
By Theorem \ref{Tdirlimit}, we may assume that $\Gamma$ is finitely
generated.  By assumption there is a short exact sequence
\[
1 \longrightarrow T \longrightarrow \Gamma
\overset{q}{\longrightarrow} C \longrightarrow 1
\]
where $T$ is a torsion abelian group and $C$ is virtually cyclic.
If $C$ was finite, then so would $\Gamma$ and the result is obvious.
Consequently we may assume that $C$ is infinite, and then there is a
short exact sequence
\[
1 \longrightarrow  Z \longrightarrow C \longrightarrow F
\longrightarrow 1
\]
where $Z$ is infinite cyclic and $F$ is finite.  The conjugation
action of $\Gamma$ on $T$ makes $T$ into
a finitely generated $\mathbb
{Z}C$-module.  Since $T$ is a torsion abelian group, it must have
finite exponent.

First we consider the case when $T$ has prime exponent $p$.  Then
$T$ is a finitely generated module over the principal ideal domain
$\mathbb {Z}_p Z$ and hence is a direct sum of finitely many cyclic
submodules.  Each cyclic module is either finite or free.  Let $T_0$
indicate the sum of the free submodules.  Then $T_0$ is a free
$\mathbb {Z}_pZ$-module of rank $n$ for some integer $n$ and
$T_0 \lhd q^{-1} (Z)$.  Choose $z_0 \in q^{-1}(Z)$ such that
$\langle q(z_0) \rangle = Z$, and then set $Z_0 = \langle z_0
\rangle$ and $\Gamma_0 = T_0 Z_0$.  Then $Z_0 \cong \mathbb {Z}$ and
$T_0$ is a free $\mathbb {Z}_pZ_0$-module of rank $n$.  Therefore
$\Gamma_0 \cong \mathbb {Z}_p^n \wr \mathbb {Z}$.
Also $\Gamma_0$ is
a subgroup of finite index in $\Gamma$, hence $\Gamma$ is isomorphic
to a subgroup of $\Gamma_0 \wr F_0$ for some finite group $F_0$ by
Lemma \ref{Lalgebraiclemma}.  Therefore $\Gamma$ is
isomorphic to a subgroup of
\[
(\mathbb {Z}_p^n \wr \mathbb {Z}) \wr F_0.
\]
The result now follows from Lemmas \ref{Lwreathp} and
\ref{Lficsubgroup}.

In general we proceed by induction on the exponent $\sigma$ of $T$.
If $\sigma = 1$ the result is obvious.  Otherwise we
may write $\sigma = p\tau$ for some prime $p$ and some integer $\tau
< \sigma$.  Note that $pT$ is a normal subgroup of $\Gamma$.
There are short exact sequences
\begin{gather*}
1 \longrightarrow pT \longrightarrow \Gamma
\overset{\phi}{\longrightarrow} \Gamma_1
\longrightarrow 1 \\
1 \longrightarrow T/pT \longrightarrow \Gamma_1
\longrightarrow C \longrightarrow 1
\end{gather*}
where $\Gamma_1 = \Gamma/pT$.  From the above, we see that
$\Gamma_1$ satisfies the FIC because the exponent of $T/pT$ is $p$.
Therefore by Lemma \ref{Lficinduction}
we need to prove the FIC holds for $\phi^{-1}(S)$ for every
virtually cyclic subgroup $S$ of $\Gamma_1$.  But $pT$ has exponent
$\tau$ and since $\tau < \sigma$, our inductive assumption implies
that $\phi^{-1}(S)$ satisfies the FIC, as required.
\end{proof}

\begin{Lem} \label{Lwreath}
Let $n$ be a non-negative integer and let $F$ be a finite group.
Then $(\mathbb {Z}^n \wr \mathbb {Z}) \wr F$ satisfies the FIC.
\end{Lem}
\begin{proof}
Fix an abelian group $A$ and a positive integer $i$.  Let $G_i =
A^{i+1}*_{A^i}$, the HNN extension with respect to the following two
embeddings of $A^i$ in $A^{i+1}$:
\[
(a_1, \dots, a_i) \mapsto (a_1, \dots, a_i,0)
\quad\text{and}\quad
(a_1, \dots, a_i) \mapsto (0, a_1, \dots, a_i).
\]
Then $A \wr \mathbb {Z} = \varinjlim G_i$ and hence
\[
(A \wr \mathbb {Z}) \wr F = \varinjlim (G_i \wr F).
\]
Consequently by Theorem \ref{Tdirlimit}, it suffices to show that
each $G_i \wr F$ satisfies the FIC when $A = \mathbb {Z}^n$.  For
the remainder of this argument we fix $A = \mathbb {Z}^n$.  Let
$\overline {G_i} = A^{i+1} \rtimes_{\alpha} \mathbb {Z}$, where
$\alpha (a_1, \dots , a_{i+1}) = (a_{i+1}, a_1, \dots, a_i)$.
There is a natural epimorphism $\phi \colon G_i \twoheadrightarrow
\overline {G_i}$.

Viewing $G_i$ as an HNN extension, we have an associated tree $T$,
and $G_i$ acts on this tree \cite[I.3.4]{dicksdunwoody}.  The
vertices of $T$ are the left cosets $gA^{i+1}$, and the edges are
the left cosets $gA^i$ for $g \in G_i$.  Also the edge $gA^i$ joins
$gA^{i+1}$ to $hA^{i+1}$, where $hA^{i+1}$ is the unique left coset
such that $gA^{i+1} \cap hA^{i+1} = gA^i$.  The stabilizer of the
vertex $gA^{i+1}$ is the subgroup $gA^{i+1}g^{-1}$, and the
stabilizer of the edge $gA^i$ is the subgroup $gA^i g^{-1}$.
Also if $e$ is an edge with stabilizer $E$ and $v$ is a vertex at
one of the ends of $e$ with stabilizer $V$, then
\begin{equation} \label{Estabilizer}
V/E \cong A.
\end{equation}

Subgroups of $G_i$ will act on this tree.  Now $\ker \phi$ acts on
$T$ and $\ker \phi \cap A^{i+1} = 1$.  Therefore $\ker \phi$ acts
freely on $T$ and we deduce from \cite[I.4.1]{dicksdunwoody} that
$\ker \phi$ is a free group.

Next consider a virtually cyclic subgroup $C$ of $\overline {G_i}$;
of course such a subgroup will be either 1 or infinite cyclic.
Then $\phi^{-1} (C)$
will act on $T$.  Since $\phi^{-1}(C)/\ker\phi \cong
1$ or $\mathbb {Z}$, we
see that the stabilizer of each vertex is either 1 or $\mathbb {Z}$,
and also that every subgroup of $\phi^{-1}(C)$ is free or free by
infinite cyclic.

Suppose the stabilizer in
$\phi^{-1}(C)$ of an edge $e$ is not 1.  Then
it must be an infinite cyclic group $C_0$.
Using \eqref{Estabilizer}, we now see that the stabilizer
of the two vertices of the ends of $e$ is also $C_0$.

We can now apply the structure theorem for groups acting on trees
\cite[I.4.1]{dicksdunwoody}.  We see that $\phi^{-1}(C)$ is the
fundamental group of a graph of groups
$\mathcal {G}$ in which all the edge groups
and vertex groups are either 1 or an infinite cyclic group.
Furthermore if an edge group is an infinite cyclic group, then the
monomorphisms into the vertex groups at the end of the edge are
onto.

The fundamental group of $\mathcal {G}$ is the ascending union of
fundamental groups of graphs of groups whose underlying graphs are
finite connected subgraphs of the underlying graph of $\mathcal
{G}$.  Let us consider one such graph of groups $\mathcal {G}_0$, so
the underlying graph of $\mathcal {G}_0$ is finite and connected.

Following \cite[p.~193]{scottwall}, we say that an
edge $e$ of $\mathcal {G}_0$ is trivial if the two ends of $e$ are
distinct vertices, and the monomorphism of the edge group into one
of the vertex groups is an isomorphism.
We can simplify $\mathcal {G}_0$ to a \emph{minimal} graph
of groups $\mathcal {G}_1$ by eliminating trivial edges
\cite[p.~193]{scottwall}.  Whenever we have an edge with
distinct vertices (i.e.\ the edge is not a loop) such that the edge
group maps isomorphically onto
one of the two vertex groups, we may contract
that edge to a point.  Then all the vertex groups of
$\mathcal {G}_1$ are 1 or infinite cyclic,
and an edge group is infinite cyclic if and
only if it is a loop and the corresponding monomorphisms into the
vertex group are isomorphisms.  We conclude that $\mathcal
{G}_1$ looks like
\bigskip

{\large
\setlength{\unitlength}{1mm}
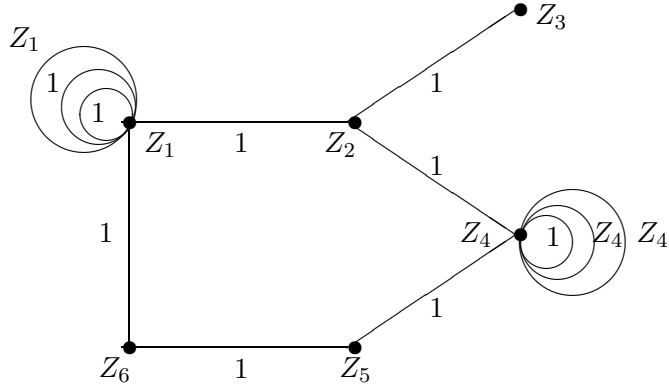
\begin{picture}(20,50)(-20,-3)
\put(0,0){$\boldsymbol \bullet$}
\put(-3,-3){$Z_6$}
\put(0,1){\line(1,0){30}}
\put(30,0){$\boldsymbol \bullet$}
\put(29,-3){$Z_5$}
\put(27,27){$Z_2$}
\put(3,27){$Z_1$}
\put(0,30){$\boldsymbol \bullet$}
\put(30,30){$\boldsymbol \bullet$}
\put(1,1){\line(0,1){30}}
\put(0,31){\line(1,0){30}}
\put(52,15){$\boldsymbol \bullet$}
\put(30,1){\line(3,2){22.5}}
\put(30,31){\line(3,-2){22.5}}
\put(30,31){\line(3,2){22.5}}
\put(52,45){$\boldsymbol \bullet$}
\put(55,44){$Z_3$}
\put(45,15){$Z_4$}
\put(-2,32){\circle{7}}
\put(56.5,15){\circle{7}}
\put(58,15){\circle{10}}
\put(60,15){\circle{20}}
\put(56.5,14.5){1}
\put(62.5,15){$Z_4$}
\put(68.5,15){$Z_4$}
\put(-4,31){1}
\put(-5,34){\circle{20}}
\put(-3,33){\circle{10}}
\put(-10,35){1}
\put(-15,41){$Z_1$}
\put(15,-3){1}
\put(15,27){1}
\put(41,35){1}
\put(41,24){1}
\put(41,5){1}
\put(-3,15){1}
\end{picture}}
\bigskip

\noindent
where $Z_i \cong \mathbb {Z}$ for all $i$.  Now
remove as many edges with trivial associated group from $\mathcal
{G}_1$ as possible so that the graph remains connected, to obtain a
new graph of groups $\mathcal {G}_2$ which has the same vertices as
$\mathcal {G}_1$.  Then the underlying graph of
$\mathcal {G}_2$ is a tree modulo the loops,
the loops have associated group $\mathbb {Z}$, and the
fundamental group of $\mathcal {G}_1$ is the free product of the
fundamental group of $\mathcal{G}_2$ and a free group.
Thus $\mathcal {G}_2$ looks like

{\large
\setlength{\unitlength}{1mm}
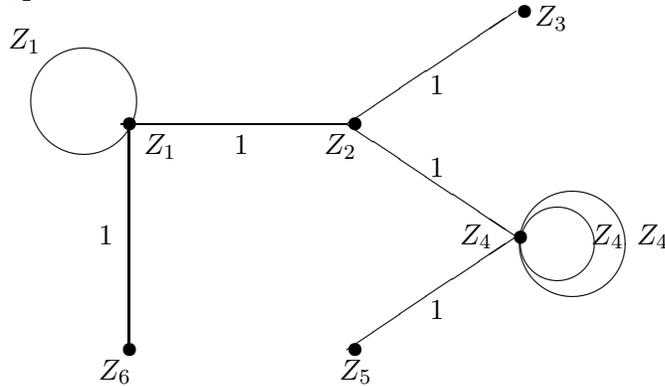
\begin{picture}(20,50)(-20,-3)
\put(0,0){$\boldsymbol \bullet$}
\put(-3,-3){$Z_6$}
\put(30,0){$\boldsymbol \bullet$}
\put(29,-3){$Z_5$}
\put(27,27){$Z_2$}
\put(3,27){$Z_1$}
\put(0,30){$\boldsymbol \bullet$}
\put(30,30){$\boldsymbol \bullet$}
\put(0,31){\line(1,0){30}}
\put(52,15){$\boldsymbol \bullet$}
\put(30,1){\line(3,2){22.5}}
\put(30,31){\line(3,-2){22.5}}
\put(30,31){\line(3,2){22.5}}
\put(52.5,45){$\boldsymbol \bullet$}
\put(55,44){$Z_3$}
\put(45,15){$Z_4$}
\put(58,15){\circle{10}}
\put(60,15){\circle{20}}
\put(62.5,15){$Z_4$}
\put(68.5,15){$Z_4$}
\put(-5,34){\circle{20}}
\put(-15,41){$Z_1$}
\put(15,27){1}
\put(41,35){1}
\put(41,24){1}
\put(41,5){1}
\put(-3,15){1}
\put(1,1){\line(0,1){30}}
\end{picture}}
\bigskip

Corresponding to each loop, we have an automorphism of the
corresponding vertex group $Z_i$.  Since $Z_i \cong \mathbb {Z}$,
this automorphism must be either the identity or inversion.  However
if it was inversion, then a nonidentity
element of the vertex group would be conjugate to its inverse in
$\phi^{-1}(C)$.  But
$\ker \phi$ intersects all vertex groups trivially and $C$ is
cyclic.  Certainly no nonidentity
element in a cyclic group can be conjugate to its
inverse, so we now have a contradiction.
Therefore each automorphism of the corresponding vertex group must
be the identity.

Now let
\[
B = \langle x,t_1, \dots, t_l \mid t_1xt_1^{-1} = x, \dots,
t_lxt_l^{-1} = x \rangle
\]
where $l$ is the maximum number of loops at a vertex.  Then the
fundamental group of $\mathcal {G}_2$ is isomorphic to a subgroup of
$B * \dots *B$, a finite free product of $B$'s, and it follows that
the fundamental group of $\mathcal {G}_1$ is also isomorphic to a
subgroup of $B *\dots *B$.

Clearly $\langle x \rangle$ is a central subgroup of $B$
and $B/\langle x\rangle$ is a
free group of rank $l$; in particular $B$ is torsion-free.
Thus $B$ is isomorphic to a direct product of a finitely generated
free group and an infinite cyclic group.  Using
Lemma \ref{Lalgebraiclemma}, we may embed $B$ in $H \wr
\mathbb {Z}_2$ where $H$ is a finitely generated
free group.

Summing up, we have shown that if $C$ is an arbitrary virtually
cyclic subgroup of $\overline{G_i}$, then
a finitely generated subgroup of $\phi^{-1}(C)$ is
isomorphic to a subgroup of $B* \dots *B$ (a finite free product of
$B$'s), and that $B$ is isomorphic to a subgroup of $H \wr
\mathbb {Z}_2$.

Now let $\psi \colon G_i \wr F \twoheadrightarrow \overline {G_i}
\wr F$ denote the epimorphism canonically determined by $\phi$.
Since $\overline{G_i} \wr F$ is virtually abelian, it satisfies the
FIC by Lemma \ref{Lficvirabelian}, so by
Lemma \ref{Lficinduction} it
will be sufficient to show that $\psi^{-1}(D)$ satisfies the FIC
whenever $D$ is a virtually cyclic subgroup of $\overline{G_i} \wr
F$.  Let $G_0$
denote the base group of $\overline {G_i} \wr F$ and let $f = |F|$.
Then $\psi^{-1}(D)/\psi^{-1}(D \cap G_0)$ is isomorphic
to a subgroup of $F$, and $\psi^{-1}(D \cap
G_0)$ is isomorphic to a subgroup of
\[
\psi^{-1} (C_1) \times \dots \times \psi^{-1} (C_f)
\]
for some virtually cyclic subgroups $C_1, \dots,
C_f$ of $\overline{G_i}$.
By two applications of Lemma \ref{Lalgebraiclemma}, we
see that any finitely generated subgroup of
$\psi^{-1}(D)$ is isomorphic to a subgroup of
$((B* \dots *B) \wr F) \wr F$.  Using Lemma \ref{Lwreathwreath}, we
deduce that any finitely generated subgroup of
$\psi^{-1}(D)$ is isomorphic to a subgroup of
\[
E := (B * \dots *B ) \wr_Y F_1
\]
for some finite set $Y$ and some group $F_1$ acting faithfully on
$Y$ (where there are finitely many $B$'s in the free product).
By Lemma \ref{Lficsubgroup} and Theorem \ref{Tdirlimit}
it will suffice to prove that $E$ satisfies the FIC.
Now we have an epimorphism $\pi \colon B * \dots *B
\twoheadrightarrow B$ determined by the identity map on each factor
of the free product, and this induces an epimorphism
\[
\theta \colon E \twoheadrightarrow B \wr_Y F_1
\leqslant (H \wr \mathbb {Z}_2) \wr_Y F_1.
\]
Applying Lemmas 
\ref{Lficfree} and \ref{Lwreathwreath}, we see that
$(H \wr \mathbb {Z}_2) \wr_Y F_1$ satisfies the FIC,
consequently by Lemmas \ref{Lficinduction} and \ref{Lficsubgroup},
it will be sufficient to
prove that $\theta^{-1} (S)$ satisfies the FIC
for every virtually cyclic subgroup
$S$ of $B \wr_Y F_1$.  Let $e = |Y|$ and let $B_0$
denote the base group of $B\wr_Y F_1$.
Then $B_0 \cong B^e$ and $\theta^{-1}(S \cap B_0)$
is a normal subgroup of index at most $|F_1|$ in $\theta^{-1}(S)$.
Since $B$ is torsion-free, $S \cap B_0$ is 1 or infinite cyclic.

Next consider $B* \dots *B$ acting on the standard tree associated
to this free product \cite[I.3.4]{dicksdunwoody}.
The stabilizers of
the vertices will be of the form $gBg^{-1}$ for $g \in G$, and the
stabilizers of the edges will be 1.  Now $\ker \pi$ intersects all
the stabilizers of the vertices trivially, consequently
if $D_1$ is an infinite cyclic subgroup of $B$, then  $\pi^{-1}
(D_1)$ intersects each vertex stabilizer in at most an infinite
cyclic group.
Using the structure theorem for groups acting on
trees \cite[I.4.1]{dicksdunwoody}, we see that $\pi^{-1}(D_1)$
is a free group and hence $\theta^{-1}(S \cap B_0)$ is isomorphic
to a subgroup of a finite direct product of $m$ copies
of a finitely generated free group $H_0$
for some positive integer $m$, i.e.\
a subgroup of $H_0^m$.  Let $F_2$ be a group of order $m$.
Using Lemma \ref{Lalgebraiclemma}, we deduce that
$\theta^{-1}(S)$ is isomorphic to a subgroup of $(H_0 \wr F_2)
\wr F_1$.  The result
now follows from Lemmas \ref{Lficsubgroup}, \ref{Lficfree} and
\ref{Lwreathwreath}.
\end{proof}

Using Lemma \ref{Lwreath}, we obtain the following result; it is
during the proof of this result where the
assumption that nearly crystallographic groups satisfy the FIC is
used.  It is the only place in this paper where this assumption is
used.
\begin{Cor} \label{Cwreath}
Assume that every nearly crystallographic group
satisfies the FIC.  Then
every torsion-free abelian by virtually cyclic group satisfies
the FIC.
\end{Cor}
\begin{proof}
Let $\Gamma$ be an extension of a virtually cyclic group $C$ by a
normal torsion-free abelian subgroup $A$.  Then we have an exact
sequence
\begin{equation} \label{Ewreath}
1 \longrightarrow A \longrightarrow \Gamma
\overset{\phi}{\longrightarrow} C \longrightarrow 1.
\end{equation}
Since $A$ is torsion-free, it embeds in $A \otimes
\mathbb {Q}$, and this embedding induces a map $H^2(C,A) \to H^2(C,
A\otimes \mathbb {Q})$.  It follows that the extension
\eqref{Ewreath} embeds in the extension
\[
1 \longrightarrow A \otimes \mathbb {Q} \longrightarrow
\hat{\Gamma} \longrightarrow C \longrightarrow 1,
\]
where $\hat{\Gamma}$ is a group containing $\Gamma$ as a subgroup.
Hence by Lemma \ref{Lficsubgroup}, we may as well assume that $A$ is
a $\mathbb {Q}$-vector space.  But then $H^2(C,A) = 0$ and so
sequence \eqref{Ewreath} splits.  Therefore we may write $\Gamma$ in
the form $A \rtimes C$.  Using Theorem \ref{Tdirlimit}, we
may assume that $A$ is finitely generated as a $\mathbb
{Q}C$-module.

Now if $C$ is finite, then $\Gamma$ is virtually abelian and the
result follows from Lemma \ref{Lficvirabelian}.  Therefore we may
assume that we have an exact sequence
\[
1 \longrightarrow Z \longrightarrow C \longrightarrow
F \longrightarrow 1
\]
where $F$ is finite and $Z \cong \mathbb {Z}$.
Then $A$ is a finitely generated $\mathbb {Q}Z$-module.  Since
$\mathbb {Q}Z$ is a principal ideal domain, $A$ is the direct sum of
a finite number of cyclic submodules.  There are two types: free and
finite $\mathbb {Q}$-dimensional.  Let $B$ denote the sum of all
the finite dimensional ones, and let $m = \dim_{\mathbb {Q}}B$.
Then $A/B$ is a direct sum of $n$ free $\mathbb {Q}Z$-modules
for some integer $n$.  Note
that $B$ is a normal subgroup of $\Gamma$, so the conjugation action
of $\Gamma$ on $B$ makes $B$ into a $\mathbb {Q}C$-module.

We shall use induction first on $n$, and then on $m$.  If $m=n=0$,
the result is obvious because then $\Gamma$ is virtually cyclic, so
assume that $n=0$ and $m > 0$.  Let $B_0$ be an irreducible
$\mathbb {Q}C$-module of $B$.
By induction $\Gamma/B_0$ satisfies the
FIC, so by Lemma \ref{Lficinduction} it will be sufficient to verify
the FIC for $\Gamma_0$ where $\Gamma_0/B_0$ is a
virtually cyclic subgroup of $\Gamma/B_0$.  Then we have an exact
sequence
\[
1 \longrightarrow B_0 \longrightarrow \Gamma_0
\longrightarrow \Gamma_0/B_0 \longrightarrow 1
\]
and by the same argument as for sequence \eqref{Ewreath}, this
sequence splits.  Using Theorem
\ref{Tdirlimit} and the assumption that nearly crystallographic
groups satisfy the FIC, we see that the FIC is true for $\Gamma_0$
(this is the only place that we need the assumption that nearly
crystallographic groups satisfy the FIC).

Therefore we may assume that $n>0$.  Next we establish the FIC for
$\Gamma/B$.  We have exact sequences
\begin{gather*}
1 \longrightarrow A/B \longrightarrow \Gamma/B \longrightarrow C
\longrightarrow 1 \\
1 \longrightarrow \phi^{-1}(Z)/B \longrightarrow \Gamma/B
\longrightarrow F \longrightarrow 1,
\end{gather*}
so $\Gamma/B$ is isomorphic to a subgroup of $(\phi^{-1}(Z)/B)
\wr F$ by Lemma \ref{Lalgebraiclemma}.
Using Lemma \ref{Lficsubgroup}, it will be sufficient to show that
$(\phi^{-1}(Z)/B) \wr F$ satisfies the FIC.  Now
$\phi^{-1}(Z)/B \cong \mathbb {Q}^n \wr \mathbb {Z}$ and $\mathbb
{Q}$ is a direct limit of subgroups isomorphic to $\mathbb {Z}$.
Therefore
\[
(\phi^{-1}(Z)/B) \wr F \cong \varinjlim (\mathbb {Z}^n \wr \mathbb
{Z}) \wr F.
\]
By Lemma \ref{Lwreath} and Theorem \ref{Tdirlimit},
the right hand side of
the above isomorphism satisfies the FIC, so we have established that
$\Gamma/B$ satisfies the FIC.

Therefore by Lemma \ref{Lficinduction}
it remains to prove the FIC when $\Gamma/B$ is virtually
cyclic.  But then we are back in the case $n=0$ and we are finished.
\end{proof}

Recall that an $n$-step solvable group means a group
$G$ with $G^{(n)}
= 1$, where $G^{(n)}$ denotes the $n$th term of the derived series,
so $G^{(0)} = G$ and $G^{(1)} = G'$.
Thus a 1-step solvable group is the same as an abelian group.
Our next step is to establish Theorem \ref{T2} for the case
of virtually two-step solvable groups.

\begin{Lem} \label{L2step}
Let $\Gamma$ be a virtually two-step solvable group and assume that
the FIC is true for nearly crystallographic groups.  Then the FIC is
true for $\Gamma$.
\end{Lem}
\begin{proof}
By assumption there is a short exact sequence $1 \to \pi \to \Gamma
\to F \to 1$ where $\pi$ is 2-step solvable
and $F$ is finite.  Since
$[\pi,\pi]$ is normal in $\Gamma$, we have short exact sequences
\begin{gather*}
1 \longrightarrow [\pi, \pi] \longrightarrow \Gamma
\overset{p}{\longrightarrow} \Gamma_0 \longrightarrow 1 \\
1 \longrightarrow A \longrightarrow \Gamma_0 \longrightarrow
F \longrightarrow 1
\end{gather*}
where $A$ is the abelian group $\pi/[\pi,\pi]$.
Since $\Gamma_0$ is virtually abelian, it satisfies the FIC by Lemma
\ref{Lficvirabelian}.  Now let $S$ be a virtually cyclic subgroup of
$\Gamma_0$ and consider the exact sequence
\[
1 \longrightarrow [\pi,\pi] \longrightarrow p^{-1}(S)
\longrightarrow S \longrightarrow 1.
\]
By Lemma \ref{Lficinduction},
it will be sufficient to show that $p^{-1}(S)$ satisfies the FIC.
Let $T$ denote the torsion subgroup of $[\pi,\pi]$, which is a well
defined abelian normal subgroup of $p^{-1}(S)$ because $[\pi,\pi]$
is an abelian normal subgroup of $p^{-1}(S)$.
Then $T \lhd
p^{-1}(S)$ and we have a short exact sequence
\[
1 \longrightarrow T \longrightarrow p^{-1}(S) \longrightarrow
p^{-1}(S)/T \longrightarrow 1.
\]
Note that $p^{-1}(S)/T$ is torsion-free abelian by virtually cyclic,
so according to
Corollary \ref{Cwreath} it satisfies the FIC.  By Lemma
\ref{Lficinduction}, it remains to prove that
$p^{-1}(C)$ satisfies the FIC whenever $C/T$ is a virtually
cyclic subgroup of $p^{-1}(S)/T$.  But now Corollary \ref{Cwreathp}
applies and we are finished.
\end{proof}

\begin{proof}[Proof of Theorem \ref{T2}]
Let $\Gamma$ be a virtually $n$-step solvable group.  We will prove
by induction on $n$ that $\Gamma$ satisfies the FIC; the inductive
step is to assume that all virtually $(n-1)$-step solvable
groups satisfy the FIC.  By Lemma
\ref{L2step}, we may assume that $n \ge 3$.
Since $\Gamma$ is assumed to be virtually $n$-step solvable, there
exists a short exact sequence
\[
1 \longrightarrow \pi \longrightarrow \Gamma \longrightarrow F
\longrightarrow 1
\]
where $F$ is finite and $\pi$ is $n$-step solvable.  Now consider
the short exact sequence
\[
1 \longrightarrow \pi^{(2)} \longrightarrow \Gamma
\overset{p}{\longrightarrow} \Gamma_0 \longrightarrow 1
\]
where $\Gamma_0 = \Gamma/\pi^{(2)}$.  Then $\Gamma_0$ is virtually
2-step solvable and hence satisfies the FIC by
Lemma \ref{L2step} (or by induction).  Let $S
\leqslant \Gamma_0$ be a virtually cyclic subgroup.
Then $p^{-1}(S)$ is virtually
$(n-1)$-step solvable and hence satisfies the FIC by our
inductive assumption.  An application of Lemma \ref{Lficinduction}
completes the proof.
\end{proof}

\section{$L$-theory analogues} \label{SLtheory}

Given a group $\Gamma$, let $L_n(\Gamma)$ denote
$L^{-\infty}_n(\mathbb {Z}\Gamma) \otimes \mathbb {Z} [\frac{1}{2}]$
where $\mathbb {Z}\Gamma$ is equipped with the standard
involution determined by the involution $\gamma \mapsto \gamma^{-1}
\colon \Gamma \to \Gamma$.  Notice that since the
difference between the functors $L^s$,
$L^h$ and $L^{-\infty}$ is only 2-torsion, we can equally well
define $L_n(\Gamma)$ by
\[
L_n(\Gamma) = L_n^s(\mathbb {Z}\Gamma) \otimes \mathbb
{Z}\bigl[\frac{1}{2}\bigr] = L_n^h(\mathbb {Z}\Gamma) \otimes
\mathbb {Z}\bigl[\frac{1}{2}\bigr].
\]
Let $\mathcal {C}$ and $\mathcal {F}$ denote
the classes of virtually
cyclic and finite groups respectively, and let $\phi \colon A \to B$
be an epimorphism between groups.
Then $\mathfrak{a}_{\mathcal{C}}$,
$\mathfrak{a}_{\mathcal {F}}$ will denote
the Quinn assembly maps in the
FIC for $\phi$ and the spectra valued functor $\underline{L}(\,)$
relative to the classes $\mathcal {C}$ and $\mathcal {F}$,
respectively; cf.\ \cite[remark A.12]{farrelljonesiso} and Remark
\ref{Raddendum} below.  Also let $\mathfrak{a}_{\mathcal {F},
\mathcal{C}}$ denote the relative assembly
map.  It is a direct consequence of Cappell's generalization of
Wall's formula for $L_n(\pi \rtimes_{\alpha} C)$ that
$\mathfrak{a}_{\mathcal {F}}$ is an
equivalence of spectra for all $S
\in \mathcal {C}$.  Hence \cite[A.12.1]{farrelljonesiso} yields the
following result.
\begin{Lem} \label{LL1}
For any epimorphism $\phi$, $\mathfrak{a}_{\mathcal {F}}$ is an
equivalence of spectra if and only if $\mathfrak{a}_{\mathcal {C}}$
is an equivalence of spectra.
\end{Lem}

\begin{Thm} \label{TL1}
Every elementary amenable group satisfies the FIC for the functor
$\underline{L}(\,)$ relative to both classes of subgroups
$\mathcal {F}$ and $\mathcal {C}$.
\end{Thm}

Given a group $\Gamma$, let
\[
\mathcal{L}_n(\Gamma) = L^{-\infty}_n(\mathbb {Q}\Gamma) \otimes
\mathbb {Z}[\frac{1}{2}]
\]
where the involution on $\mathbb {Q}\Gamma$ is determined by
the involution $\gamma \mapsto \gamma^{-1} \colon
\Gamma \to \Gamma$.

\begin{Cor} \label{CL1}
Every elementary amenable group satisfies the FIC for the functor
$\underline{\mathcal {L}}(\,)$ relative to both classes of subgroups
$\mathcal {F}$ and $\mathcal {C}$.
\end{Cor}

Before proving Theorem \ref{TL1}, we deduce Corollary \ref{CL1} from
it.  The key fact, due to Ranicki \cite{ranicki},
is that the natural map
\[
L^{-\infty}_n(\mathbb {Z}\Gamma) \longrightarrow
\mathcal {L}^{-\infty}_n (\mathbb {Q}\Gamma)
\]
is an isomorphism after tensoring with $\mathbb {Z}[\frac{1}{2}]$.
Hence the natural map $\underline {L} (\,) \to
\underline {\mathcal {L}}(\,)$ is
an equivalence of spectra valued functors.  Thus Corollary \ref{CL1}
is clearly equivalent to Theorem \ref{TL1}.

\begin{proof}[Proof of Theorem \ref{TL1}]
We shall use a standard inductive procedure for elementary amenable
groups.  Thus let $\mathcal {A}$ denote
the class of virtually finitely
generated abelian groups, and if $\mathcal {X}$ and $\mathcal {Y}$
are classes of groups, then $G \in \mathcal {X}\mathcal {Y}$ will
mean that there exists $H \lhd G$ such that $H \in \mathcal {X}$ and
$G/H \in \mathcal {Y}$, and $G \in L\mathcal {X}$ will mean that
every finitely generated subgroup of
$G$ is contained in an $\mathcal
{X}$-group (if $\mathcal {X}$ is closed under taking subgroups, then
this is equivalent to saying that every finitely generated subgroup
of $G$ is an $\mathcal {X}$-group).  Then for each ordinal $\alpha$,
define the class of groups $\mathcal {E}_{\alpha}$ inductively as
follows:
\begin{align*}
\mathcal {E}_0 &= \mathcal {A}, \\
\mathcal {E}_{\alpha} &= \bigcup_{\beta < \alpha} \mathcal
{E}_{\beta} \quad \text{if }\alpha\text{ is a limit ordinal,}\\
\mathcal {E}_{\alpha} &= (L\mathcal {E}_{\beta})\mathcal {A}\quad
\text{if }\alpha\text{ is a successor ordinal and }\alpha= \beta+1.
\end{align*}
Then $\bigcup_{\alpha \ge 0} \mathcal {E}_{\alpha}$ is the class of
elementary amenable groups.  Furthermore if $H \lhd G$, $G/H$ is
finite, and $H \in \mathcal {E}_{\alpha}$ 
or $L\mathcal {E}_{\alpha}$, then $G \in
\mathcal {E}_{\alpha}$ or $L\mathcal {E}_{\alpha}$
respectively, cf.\ \cite[\S 3]{klm}.

Now let $\Gamma$ be the elementary amenable group for which we want
to prove the theorem.  If $\Gamma \in \mathcal{E}_0$, then $\Gamma$
is a virtually finitely generated abelian group, and the result for
such groups is true by \cite[theorem 2.1 and remark
2.1.3]{farrelljonesiso}.  Therefore
we may assume that $\Gamma \notin \mathcal {E}_0$.
Let $\alpha$ be the least ordinal such that $\Gamma
\in \mathcal {E}_{\alpha}$.  Now $\alpha$ cannot be a limit ordinal,
hence it must be a successor ordinal and we may write $\alpha =
\beta +1$ for some ordinal $\beta$.
Then there is an epimorphism $\phi \colon \Gamma \to D$ where
$\ker \phi \in L\mathcal{E}_{\beta}$ and $D \in \mathcal {A}$.
At this point we need the following analogue of \cite[proposition
2.2]{farrelljonesiso}.

\begin{Lem} \label{LL2}
Let $\psi \colon A \twoheadrightarrow B$ be a group epimorphism.
Suppose the FIC relative to the functor $\underline{L}(\,)$ and the
class $\mathcal {F}$ is true for $B$ and $\psi^{-1}(G)$
for all finite subgroups $G$ of $B$.  Then the FIC is true for
$A$.
\end{Lem}
We omit the proof of Lemma \ref{LL2} because it is virtually
identical to the one given
for \cite[proposition 2.2]{farrelljonesiso}; cf.\
\cite[remark A.12]{farrelljonesiso}.

We now apply Lemma \ref{LL2} (with $\psi = \phi$, $B=D$ and $A =
\Gamma$) to complete the inductive step in the proof of Theorem
\ref{TL1}.  Note that $D$ satisfies the FIC because $D \in \mathcal
{A}$.  Now set $\pi = \ker \phi$ and
$\tilde {\pi} = \phi^{-1}(G)$ where $G$ is a finite
subgroup of $D$.  Since $\pi \in L\mathcal {E}_{\beta}$ and
$\tilde{\pi}/\pi$ is finite, we see that $\tilde{\pi} \in L\mathcal
{E}_{\beta}$ and thus $\phi^{-1}(G)$ satisfies the FIC by Theorem
\ref{Tdirlimit}.  We hence conclude from Lemma \ref{LL2} that
$\Gamma$ satisfies the FIC.
\end{proof}

\section{Application to infinite dimensional $L$-groups}
\label{SLtheoryapplication}

Wall \cite{Wall76} showed for
every finite group $\Gamma$, that the abelian groups $L_i^s(\mathbb
{Z}\Gamma)$ and $L_i^h(\mathbb {Z}\Gamma)$ are finitely generated for
all indices $i$.  However Cappell \cite{Cappell74}
using his Unil-functor
showed that for the infinite dihedral group $\mathcal {D}$, neither
$L_2^s(\mathbb {Z} \mathcal {D})$ nor $L_2^h
(\mathbb {Z} \mathcal {D})$ are finitely generated.
But the following query seemed possibly true.
\begin{Ques} \label{QBaumslag}
Is $\dim_{\mathbb {Q}} L_i (\mathbb {Z}\Gamma) \otimes \mathbb {Q} <
\infty$ for every finitely presented group $\Gamma$ and every index
$i$?
\end{Ques}
This was answered negatively by Weinberger in \cite{weinberger}; we
would like to thank Andrew Ranicki for pointing out this reference to
us.  However we shall show that our results give further negative
examples for Question \ref{QBaumslag}.

The method is to provide examples of finitely
presented groups which
have non-vanishing homology in arbitrary high dimension.
For the rest of this section, let $G$ be the group of
\cite[Example 1]{baumslagdyer}, i.e.
\[
\langle a,s,t \mid s^t = s,\ [a^t,a] = 1,\ a^s = a a^t \rangle.
\]
Here $s^t = tst^{-1}$ and $[a^t,a] = 1$ means that $a^t$ commutes
with $a$.  Then \cite{baumslagdyer} shows that $G$ has the
following properties:

\begin{enumerate} [\normalfont (i)]
\item
$G$ is generated by three elements.

\item
$G$ is finitely presented torsion-free metabelian.

\item \label{infinite}
$H_n(G,\mathbb {Z})$ is free abelian of infinite rank for all
$n \ge 3$.

\end{enumerate}

\begin{proof}[Proof that Question \ref{QBaumslag} has a negative
answer]
By Theorem \ref{TL1}, for any
elementary amenable group $\Gamma$
\[
\mathbb {H}_i(B \Gamma, \underline{\mathbb {L}}(\mathbb {Z})
\otimes \mathbb {Q}) \cong
L_i(\mathbb {Z \Gamma}) \otimes \mathbb {Q}
\]
for all $i \in \mathbb {Z}$.  But for any space $X$
\[
\mathbb {H}_i(X, \underline{\mathbb {L}}
(\mathbb {Z}) \otimes \mathbb {Q}) \cong
\bigoplus_{j= -\infty}^{\infty} H_{i + 4j}(X, \mathbb {Q}),
\]
where for convenience we interpret $H_r(Y, \mathbb {Q}) = 0$ for
$r<0$ and $Y$ is any space or group. 
Concatenating these two isomorphisms yields
\[
L_i (\mathbb {Z} \Gamma) \otimes \mathbb {Q} \cong
\bigoplus_{j= -\infty}^{\infty} H_{i+4j}(\Gamma, \mathbb {Q}).
\]
Of course \eqref{infinite} above
tells us that $H_n(G,\mathbb {Q})$ is
an infinite dimensional $\mathbb {Q}$-vector space for all $n \ge
3$, so in the case $\Gamma = G$, we see that $\dim_{\mathbb{Q}}
L_i(\mathbb {Z}\Gamma) \otimes \mathbb {Q} = \infty$ for all $i$.
Furthermore $\Gamma$ is a three generated finitely presented
torsion-free metabelian group, so we have the required example.
\end{proof}

We end this section
by noting that a more sophisticated negative answer to
Question \ref{QBaumslag} is furnished by letting $\Gamma$ be
Thompson's group $F$ of orientation preserving piecewise linear
homeomorphisms of $\mathbb {R}$, as described in
\cite{BrownGeoghegan84}.  Dan Farley \cite{Farley00}
has recently shown that
Novikov's conjecture is true for $F$; consequently
\[
\bigoplus_{j=-\infty}^{\infty} H_{i+4j}(BF, \mathbb {Q}) \subseteq
L_i(\mathbb {Z}F) \otimes \mathbb {Q}.
\]
Also Brown and Geoghegan have calculated the homology of $F$; in
particular it follows from \cite[Theorem 7.1]{BrownGeoghegan84} that
\[
H_j(BF, \mathbb {Q}) \cong \mathbb {Q} \oplus \mathbb {Q}
\quad\text{for all }j>0.
\]

\section{Appendix} \label{Sappendix}

We start by recalling the previous statement of the Fibered
Isomorphism Conjecture (FIC) made in \cite{farrelljonesiso}.  Let
$\mathcal {S}$ be a homotopy invariant (covariant) functor from the
category of topological spaces to spectra.  Important examples of
such functors are the stable topological pseudo-isotopy functor
$\mathcal {P(\,)}$, the algebraic $K$-theory functor $\mathcal
{K}(\,)$, and the $L$-theory functor
$\mathcal{L}^{-\infty} (\,)$; cf.\
\cite[\S 1]{farrelljonesiso}.  Let $\mathcal {M}$ denote the
category
of continuous surjective
maps; i.e.\ an object in $\mathcal {M}$ is a continuous
map $p \colon E \to B$ between topological spaces $E$ and $B$, while
a morphism from $p_1 \colon E_1 \to B_1$ to $p_2 \colon E_2 \to B_2$
is a pair of continuous maps $f\colon E_1 \to E_2$, $g \colon B_1
\to B_2$ making the following diagram a commutative square of maps:
\[
\begin{CD}
E_1 @>f>> E_2 \\
@Vp_1VV @VVp_2V \\
B_1 @>g>> B_2.
\end{CD}
\]
Quinn \cite[appendix]{quinninv} constructed a functor from $\mathcal
{M}$ to the category of $\Omega$-spectra which associates to the map
$p$ the spectrum $\mathbb {H}(B; \mathcal {S}(p))$ in such a way
that
\[
\mathbb {H}(B; \mathcal {S}(p)) = \mathcal {S}(E)
\]
in the
special case that $B$ is a single point $*$.  Furthermore the map of
spectra $\mathfrak{a} \colon \mathbb {H}(B; \mathcal {S}(p)) \to
\mathcal {S}(E)$ functorially associated to the commutative square
\[
\begin{CD}
E @>\id>> E \\
@VpVV @VVV \\
B @>>> *
\end{CD}
\]
is called the (Quinn) assembly map.

Let $\mathcal {E}$ be a $\Gamma$-space which is universal for the
class of all virtually cyclic subgroups of $\Gamma$ and let
$\mathcal {B}$ denote its orbit space $\mathcal {E}/\Gamma$; cf.\
\cite[appendix]{farrelljonesiso}.  Let $Y$ be an arbitrary free and
properly discontinuous $\Gamma$-space, and let $p \colon
Y \times_{\Gamma}  \mathcal{E}
\to \mathcal {E}/\Gamma = \mathcal {B}$ be the
continuous map determined by projection onto the second factor of
$Y \times \mathcal {E}$.  The FIC for $\mathcal {S}$ and $\Gamma$ is
the assertion that the assembly map
\[
\mathfrak{a} \colon \mathbb {H}(\mathcal {B}; \mathcal {S}(p))
\longrightarrow \mathcal {S}(Y \times_{\Gamma} \mathcal{E}) =
\mathcal {S} (Y/\Gamma)
\]
is a (weak) equivalence of spectra.  This conjecture was made in
\cite[\S 1.7]{farrelljonesiso} for the three functors
$\mathcal {S} = \mathcal {P}$, $\mathcal {K}$ and
$\mathcal{L}^{-\infty}$, and every group
$\Gamma$.  Throughout this appendix, $\mathcal {S}$ will denote any
of these three spectra valued functors
$\mathcal {P}$, $\mathcal {K}$,
$\mathcal{L}^{-\infty}$, or either of the two variants of $\mathcal
{L}^{-\infty}$ constructed in Section \ref{SLtheory}; namely
$\underline {L}$ and $\underline {\mathcal {L}}$.

\begin{Thm} \label{Tdirlimit}
Let $I$ be a directed set, and let $\Gamma_n$, $n \in I$,
be a directed system of groups with $\Gamma =
\varinjlim_{n\in I} \Gamma_n$; i.e.\ $\Gamma$ is the direct limit
of the $\Gamma_n$.  If each group $\Gamma_n$ satisfies FIC,
then $\Gamma$ also satisfies FIC.
\end{Thm}
\begin{proof}
Let $\mathcal {E}_n, \mathcal{E}$ denote universal
$(\mathcal{C}, \Gamma_n)$, $(\mathcal {C}, \Gamma)$ spaces where
$\mathcal {C}$ is the class of all virtually cyclic groups.
Also let $\mathcal {B}_n, \mathcal{B}$ denote the orbit spaces
$\mathcal {E}_n/\Gamma_n, \mathcal {E}/\Gamma$ respectively.
Let $E_n,E$ denote universal $\Gamma_n, \Gamma$
spaces for the class containing only the trivial group, and let
$B_n, B$ denote the orbit spaces $E_n/\Gamma_n, E/\Gamma$
respectively.  These are of course Eilenberg-Mac\,Lane spaces
$K(\Gamma_n,1)$, $K(\Gamma,1)$ respectively.  Let $\phi_{mn} \colon
\Gamma_m \to \Gamma_n$
($n \ge m$) denote the homomorphisms associated to the
directed system of groups $\Gamma_n$
and let $\phi_n \colon \Gamma_n \to \Gamma$
be the group homomorphisms making
\[
\Gamma = \varinjlim_{n\in I} \Gamma_n.
\]
Let $\bar{\phi}_n \colon B_n \to B$ and $\bar{\phi}
_{mn}$ be continuous induced maps; these can be chosen
functorially, i.e.\ so that
\[
\bar{\phi}_{ns} \circ \bar{\phi}_{mn} = \bar{\phi}_{ms}
\text{ and }
\bar{\phi}_{n} \circ \bar{\phi}_{mn} = \bar{\phi}_{m}
\]
where $s \ge n \ge m$.  Let $Y \to B$ be a fiber bundle
and let $Y_n \to B_n$ denote its pullback via the maps
$\bar{\phi}_n \colon B_n \to B$.  Let $\tilde{Y}_n \to
Y_n$ be the $\Gamma_n$ regular covering space which is the
pullback of the universal covering space $E_n \to B_n$ via
the map $Y_n \to B_n$, and let $p_n \colon \tilde{Y}_n
\times_{\Gamma_n} \mathcal {E}_n \to \mathcal {B}_n$ be the
simplicially stratified map determined by projection onto the second
factor of $\tilde {Y}_n \times \mathcal {E}_n$.  Likewise let $p
\colon \tilde {Y} \times_{\Gamma} \mathcal {E} \to \mathcal {B}$ be
the similarly determined simplicially stratified map.  There are
\cite[appendix]{quinninv} assembly maps between spectra
\begin{align*}
\mathfrak{a}_n &\colon \mathbb {H}(\mathcal {B}_n,
\mathcal {S}(p_n))
\longrightarrow \mathcal {S}(Y_n), \quad n \in I, \\
\mathfrak {a} &\colon \mathbb {H} (\mathcal {B}, \mathcal {S}(p))
\longrightarrow \mathcal {S}(Y).
\end{align*}
We have assumed that each $\mathfrak {a}_n$, $n \in I$, is an
equivalence of spectra and to prove Theorem \ref{Tdirlimit},
we must show that $\mathfrak {a}$ is also an equivalence of spectra.
It is of course sufficient to show that $\mathfrak {a}$ induces
isomorphisms on all the homotopy groups $\pi_i$; i.e.\ that
$\pi_i(\mathfrak {a}) \colon
\mathbb {H}_i (\mathcal {B}, \mathcal {S}(p))
\to \pi_i (\mathcal {S}(Y))$
is an isomorphism for each $i \in \mathbb {Z}$.
Consider the following commutative diagram:
\begin{equation} \label{EAcomlim}
\begin{CD}
\varinjlim_{n \in I}
\mathbb {H}_i(\mathcal {B}_n, \mathcal {S}(p_n))
@>>>
\mathbb {H}_i (\mathcal {B}, \mathcal {S}(p)) \\
@V\pi_i(\mathfrak {a}_n)VV
@VV\pi_i(\mathfrak {a})V \\
\varinjlim_{n \in I} \pi_i (\mathcal {S}(Y_n))
@>>> \pi_i (\mathcal {S}(Y)).
\end{CD}
\end{equation}
The left vertical arrow is an isomorphism by assumption.  We proceed
to show that both the top and bottom horizontal arrows are
epimorphisms.  Here are some useful remarks.

\begin{Rem} \label{Rfinpres}
Let $G$ be a finitely presented group.  Then any homomorphism $\psi
\colon G \to \Gamma$ lifts to some $\Gamma_n$; i.e.\ there exists a
homomorphism $\hat{\psi} \colon G \to \Gamma_n$
such that $\phi_n \circ \hat{\psi} = \psi$.
\end{Rem}

\begin{Rem} \label{Rfibration}
Let $q \colon Y \to B$ denote the projection map in the fiber
bundle.  Then $\pi_i(\mathcal S(Y)) = \varinjlim \pi_i(\mathcal
{S}(q^{-1}(K)))$ for all $i \in \mathbb {Z}$,
where $K$ varies over all finite subcomplexes of
$B$.  Moreover $\mathbb {H}_i(\mathcal {B}, \mathcal {S}(p)) =
\varinjlim \mathbb {H}_i(K, \mathcal {S}(p |_{p^{-1}(K)}))$, where
$K$ now varies over all finite subcomplexes of $\mathcal {B}$.
\end{Rem}

\begin{Rem} \label{Rfincomplex}
If $K$ is a finite complex, then $\pi_1(K)$ is finitely presented.
\end{Rem}

Let $\hat {\phi}_n$ denote the top horizontal map in the Cartesian
square
\[
\begin{CD}
Y_n @>\hat{\phi}_n>> Y \\
@VVV @VVqV \\
B_n @>\bar{\phi}_n>> B
\end{CD}
\]
and for each finite subcomplex $K$ of $B$, let $\sigma_K \colon
q^{-1}(K) \to Y$ denote the inclusion.  Then by applying Remarks
\ref{Rfinpres} and \ref{Rfincomplex} to the finite complex $K$, one
constructs a map $\tau_K \colon q^{-1} (K) \to Y_n$ (for some
$n \in I$) such that $\hat{\phi}_n \circ \tau_K \simeq \sigma_K$.
Combining this with Remark \ref{Rfibration} yields that the bottom
arrow in diagram \eqref{EAcomlim} is an epimorphism.  Consequently
$\pi_i(\mathfrak {a})$ is also an epimorphism.  To prove Theorem
\ref{Tdirlimit}, it
remains to show that $\pi_i(\mathfrak {a})$
is monic; but the argument for this is much more complicated.  The
first major step in this argument is showing that the top arrow in
diagram \eqref{EAcomlim} is an epimorphism.

For this purpose we introduce an auxiliary direct limit.  Let
$\bar{p}_n \colon \tilde{Y} \times_{\Gamma_n}
\mathcal {E}_n \to \mathcal {B}_n$ be
the simplicially stratified map determined by projection onto the
second factor of $\tilde {Y} \times \mathcal {E}_n$ (here $\Gamma_n$
acts on $\tilde {Y}$ via $\phi_n \colon \Gamma_n \to \Gamma$).  Then
we have a two step ladder of commutative diagrams
\begin{equation} \label{EAcomsq1}
\begin{CD}
\tilde{Y}_n \times_{\Gamma_n} \mathcal {E}_n @>>>
\tilde{Y} \times_{\Gamma_n} \mathcal {E}_n @>>>
\tilde{Y} \times_{\Gamma} \mathcal {E} \\
@Vp_nVV @VV\bar{p}_nV @VVpV \\
\mathcal {B}_n @>\id>> \mathcal{B}_n @>>> \mathcal {B}.
\end{CD}
\end{equation}
Let
\begin{align*}
\Phi_n &\colon \mathbb {H}_i(\mathcal {B}_n, \mathcal
{S}(p_n)) \longrightarrow \mathbb {H}_i(\mathcal {B}_n,
\mathcal {S} (\bar{p}_n)) \\
\text{and } \Psi_n &\colon \mathbb {H}_i (\mathcal {B}_n, \mathcal
{S}(\bar{p}_n)) \longrightarrow \mathbb {H}_i(\mathcal {B},
\mathcal {S}(p))
\end{align*}
be the homomorphisms which are functorially determined by the first
and second squares respectively, in the ladder \eqref{EAcomsq1} (see
again \cite[appendix]{quinninv}).  Then the direct
limit of the composites $\Psi_n \circ \Phi_n \colon \mathbb
{H}_i(\mathcal {B}_n, \mathcal {S} (p_n)) \to \mathbb {H}_i
(\mathcal {B}, \mathcal {S} (p))$ is the top arrow in
diagram \eqref{EAcomlim}.  Of course the direct limit of the family
$\Phi_n$ determines a homomorphism
\[
\Phi \colon \varinjlim_{n\in I}
\mathbb {H}_i(\mathcal {B}_n, \mathcal
{S}(p_n)) \longrightarrow \varinjlim_{n\in I}
\mathbb {H}_i (\mathcal {B}_n,
\mathcal {S} (\bar {p}_n))
\]
and likewise the direct limit of the family $\Psi_n$ gives a
homomorphism
\[
\Psi \colon \varinjlim_{n\in I}
\mathbb {H}_i (\mathcal{B}_n, \mathcal
{S}(\bar {p}_n)) \longrightarrow \mathbb {H}_i (\mathcal {B},
\mathcal {S} (p)).
\]
Then the top arrow of diagram \eqref{EAcomlim} is the composite
$\Psi \circ \Phi$.  Hence to show that the top arrow of diagram
\eqref{EAcomlim} is a surjection, it suffices to verify that both
$\Phi$ and $\Psi$ are surjections.

\begin{Rem} \label{Rdirabelian}
The maps in the directed system of abelian groups $\mathbb
{H}_i(\mathcal {B}_n, \mathcal {S}(p_n))$ are constructed
as follows where $n \ge m$.  Let $f_{mn} \colon \mathcal {E}_m \to
\mathcal {E}_n$ be a $\Gamma_m$-equivariant
map where $\Gamma_m$ acts on
$\mathcal {E}_n$ via the homomorphism $\phi_{mn} \colon \Gamma_m \to
\Gamma_n$.  Such a map exists because of
\cite[theorem A.2]{farrelljonesiso} and is unique up to a
$\Gamma_m$-equivariant homotopy.  There is also a
$\Gamma_m$-equivariant map $\tilde{\phi}_{mn} \colon \tilde{Y}_m \to
\tilde {Y}_n$ which is the top arrow in the Cartesian square
\[
\begin{CD}
\tilde {Y}_m @>\tilde{\phi}_{mn}>> \tilde{Y}_n \\
@VVV @VVV \\
E_m @>>> E_n
\end{CD}
\]
whose bottom arrow is the functorially defined
$\Gamma_m$-equivariant map which induces $\bar{\phi}_{mn} \colon B_m
\to B_n$.  Then $\tilde{\phi}_{mn} \times f_{mn}$ and $f_{mn}$
induce the top and bottom arrows, respectively,
in the following commutative square
\[
\begin{CD}
\tilde{Y}_m \times_{\Gamma_m}
\mathcal {E}_m @>>> \tilde{Y}_n \times_{\Gamma_n}
\mathcal {E}_n \\
@Vp_mVV @VVp_nV \\
\mathcal {B}_m @>>> \mathcal {B}_n.
\end{CD}
\]
This square functorially determines the desired group homomorphism
\[
\mathbb {H}_i (\mathcal {B}_m, \mathcal{S}(p_m)) \longrightarrow
\mathbb {H}_i (\mathcal {B}_n, \mathcal{S}(p_n)).
\]
Note that this homomorphism is independent of the choice of $f_{mn}$
because $\mathbb {H}_i(\ ,\ )$ is a homotopy functor.

Also let $f_n \colon \mathcal {E}_n \to \mathcal {E}$ be a
$\Gamma_n$-equivariant map and let $\tilde {\phi}_n \colon \tilde
{Y}_n \to \tilde {Y}$ be the $\Gamma_n$-equivariant map which is the
top arrow in the Cartesian square
\begin{equation} \label{EAcomsq2}
\begin{CD}
\tilde{Y}_n @>>> \tilde {Y} \\
@VVV @VVV \\
E_n @>>> E
\end{CD}
\end{equation}
whose bottom arrow is functorially determined by $\phi_n \colon
\Gamma_n \to \Gamma$.  Then $\tilde {\phi}_n \times f_n$ and $f_n$
determine for $n \in I$
the homomorphisms $\mathbb {H}_i(\mathcal {B}_n, \mathcal
{S}(p_n)) \to \mathbb {H}_i(\mathcal {B}, \mathcal
{S}(p))$ whose direct limit is the top arrow of diagram
\eqref{EAcomlim}.  Also the top arrow in the first square of ladder
\eqref{EAcomsq1} is induced by $\tilde{\phi}_n \times \id$
and the top arrow in the second square by $\id \times f_n$.
The bottom arrow in the second square is induced by $f_n$.

The homomorphisms in the directed system $\mathbb {H}_i( \mathcal
{B}_n, \mathcal {S}(\bar{p}_n))$ are determined by the
following commutative square where $n \ge m$:
\[
\begin{CD}
\tilde{Y} \times_{\Gamma_m} \mathcal {E}_m @>>>
\tilde {Y} \times_{\Gamma_n} \mathcal {E}_n \\
@V\bar{p}_mVV @VV\bar{p}_nV \\
\mathcal {B}_m @>>> \mathcal {B}_n
\end{CD}
\]
whose top arrow is induced by $\id \times f_{mn}$ and the bottom
arrow by $f_{mn}$.
\end{Rem}

We now proceed to show that $\Phi$ is a surjection.  To do this
consider another simplicially stratified map
\[
p_{mn} \colon \tilde {Y}_n \times_{\Gamma_m} \mathcal {E}_m
\longrightarrow \mathcal {B}_m
\]
determined by projection onto the second factor of $\tilde {Y}_n
\times \mathcal {E}_m$, which is defined for each pair of indices $n
\ge m$.  Note that $p_{mm} = p_m$.

\begin{Not}
If $p \colon A\to B$ is a continuous map and $\mathcal {K}$ is
a subspace of $B$, denote the restriction maps
\[
p |_{p^{-1} (\mathcal {K})} \colon p^{-1}
(\mathcal {K}) \longrightarrow \mathcal {K}
\]
also by $p$.
\end{Not}
Consider the commutative square
\begin{equation} \label{EAcomsq3}
\begin{CD}
\tilde{Y}_n \times_{\Gamma_m} \mathcal {E}_m
@>>>
\tilde{Y} \times_{\Gamma_m} \mathcal {E}_m \\
@Vp_{mn}VV @VV\bar{p}_mV \\
\mathcal {B}_m @>\id>> \mathcal {B}_m
\end{CD}
\end{equation}
whose top arrow is induced by $\tilde{\phi}_n \times \id$.  It
determines functorially a homomorphism
$\mathbb {H}_i(\mathcal {K}, \mathcal {S}(p_{mn})) \to \mathbb
{H}_i(\mathcal {K}, \mathcal{S} (\bar{p}_n))$ for each finite
subcomplex $\mathcal {K}$ of $\mathcal {B}_m$.

\begin{Lem} \label{Lcomsq}
Given a finite subcomplex $\mathcal {K}$ of $\mathcal {B}_m$, there
exists an index $n \ge m$ such that the homomorphism $\mathbb
{H}_i(\mathcal {K}, \mathcal {S}(p_{mn})) \to \mathbb {H}_i(\mathcal
{K}, \mathcal {S}(\bar{p}_n))$ determined by square
\eqref{EAcomsq3} is an isomorphism.
\end{Lem}
\begin{proof}
Restricted to each open simplex $\sigma$ in $\mathcal {K}$, both
$p_{mn}$ and $\bar{p}_{m}$ are fiber bundles and the
top arrow in square \eqref{EAcomsq3} is a bundle map.  We proceed
to show that, restricted to each such fiber, the top arrow is a
homotopy equivalence provided $n$ is chosen sufficiently large, and
then the lemma follows from this assertion.

Since there are only finitely many simplices in $\mathcal {K}$, it
suffices to prove our assertion for one such simplex $\sigma$.  Let
$\hat{\sigma}$ be a simplex in
$\mathcal {E}_m$ ``lying over $\sigma$" under
the map $\mathcal {E}_m \to \mathcal {E}_m /\Gamma_m = \mathcal
{B}_m$, and let $S$ be the subgroup of $\Gamma_m$ fixing
$\hat{\sigma}$.  Recall that $S$ is a virtually cyclic group.  Now
the fibers of $p_{mn}$ and $\bar{
p}_m$ are respectively
$\tilde {Y}_n/\phi_{mn}(S)$ and $ \tilde {Y}/\phi_m(S)$.
Furthermore, the following commutative square is induced from square
\eqref{EAcomsq2}
\begin{equation} \label{EAcomsq4}
\begin{CD}
\tilde{Y}_n/\phi_{mn}(S)
@>>>
\tilde{Y}/\phi_m(S) \\
@VVV @VVV \\
E_n/\phi_{mn}(S) @>>> E/\phi_m(S),
\end{CD}
\end{equation}
and the top arrow in this square is the restriction of the top arrow
of square \eqref{EAcomsq3}, while the bottom arrow of square
\eqref{EAcomsq4} is functorially induced by $\phi_n$.  Note
that this bottom arrow
is a map between Eilenberg-Mac\,Lane spaces $K(\phi_{mn}(S),1)$
and $K(\phi_m(S),1)$, since $E_n$ and $E$ are both contractible.
Hence the bottom arrow of square \eqref{EAcomsq4} is a homotopy
equivalence precisely when the epimorphism
\[
\phi_n \colon \phi_{mn}(S) \longrightarrow \phi_m(S)
\]
is an isomorphism.  But this happens provided $n$ is sufficiently
large since every subgroup of $S$ is finitely generated and
$\varinjlim \Gamma_n = \Gamma$.  Also when this happens
\eqref{EAcomsq4} is a Cartesian square,
since square \eqref{EAcomsq2}
was Cartesian and hence the top arrow of square \eqref{EAcomsq4} is
also a homotopy equivalence.
\end{proof}

Now a routine direct limit argument using
Lemma \ref{Lcomsq} together
with Remark \ref{Rfibration} shows that $\Phi$ is a surjection (the
details of this argument are left as an exercise).

It remains to show that $\Psi$ is a surjection.  For this purpose we
introduce some terminology.

\begin{Not}  Given a subspace $\mathcal {K}$ of $\mathcal {B}$, let
$\hat{\mathcal {K}}$ denote its inverse image under the projection
$\mathcal {E} \to \mathcal {B}$.
\end{Not}

\begin{Prop} \label{Pfinsubcomplex}
Given a finite subcomplex $\mathcal{K}$ of $\mathcal {B}$, there
exists an index $n \in I$ and a
$\Gamma$-equivariant map $\psi \colon
\hat{\mathcal{K}} \to \Gamma \times_{\Gamma_n} \mathcal {E}_n$.
\end{Prop}
Before proving this proposition, we first show how it implies that
$\Psi$ is surjective.  Define a $\Gamma$-equivariant map
\[
\bar{f}_n \colon \Gamma \times _{\Gamma_n} \mathcal {E}_n
\longrightarrow \mathcal {E}
\]
by the equation $\bar {f}_n ([\gamma, x]) = \gamma f_n(x)$ for
all $(\gamma, x) \in \Gamma\times \mathcal {E}_n$.  Here $[\gamma,
x]$ denotes the equivalence class in $\Gamma \times_{\Gamma_n}
\mathcal {E}_n$ of the pair $(\gamma, x)$.  Then the composite
$\bar{f}_n \circ \psi \colon \hat{\mathcal{K}} \to \mathcal {E}$ is
also $\Gamma$-equivariant and is hence $\Gamma$-equivariantly
homotopic to the inclusion $\hat{\mathcal {K}} \subseteq \mathcal
{E}$ by \cite[theorem A.2]{farrelljonesiso}.
(Notice that $\hat{K}$ is
cellular $\mathcal {C}$-free and that $\mathcal {E}$ is $\mathcal
{C}$-contractible where $\mathcal {C}$ is the class of all virtually
cyclic groups.)

Consider the commutative square
\begin{equation} \label{EAcomsq5}
\begin{alignat}{4}
\tilde{Y} \times&_{\Gamma} \, \hat{\mathcal {K}} &
&  -\hspace{-1.2ex} -\hspace{-1.2ex}\longrightarrow\ &
\tilde{Y} \times_{\Gamma}(\Gamma \times_{{\Gamma}_n} \mathcal{E}_n)
&=& \tilde {Y} &\times_{\Gamma_n} \mathcal {E}_n \notag \\
p\ &\Bigg |_{\displaystyle\hspace{-.97ex} \downarrow} & &&&&&
\Bigg |_{\displaystyle\hspace{-.97ex} \downarrow} \ \bar{p}_n  \\
&\mathcal {K}&
&-\hspace{-1.2ex} -\hspace{-1.2ex} -\hspace{-1.2ex}-
\hspace{-1.2ex}-\hspace{-1.2ex} -\hspace{-1.2ex} -\hspace{-1.2ex}
&-\hspace{-1.2ex} -\hspace{-1.2ex} -\hspace{-1.2ex}
-\hspace{-1.2ex} -\hspace{-1.2ex} -\hspace{-1.2ex}
-\hspace{-1.2ex} -\hspace{-1.2ex} -\hspace{-1.2ex}
-\hspace{-1.2ex} -\hspace{-1.2ex} -\hspace{-1.2ex}
-\hspace{-1.2ex} -\hspace{-1.2ex} -\hspace{-1.2ex}
&\rightarrow&
&\mathcal{B}_n
\notag
\end{alignat}
\end{equation}
whose top arrow is induced by $\id \times \psi$ and bottom arrow by
$\psi$.
Square \eqref{EAcomsq5} functorially determines a homomorphism
$\hat{\psi} \colon \mathbb {H}_i (\mathcal {K}, \mathcal
{S}(p)) \to \mathbb {H}_i(\mathcal {B}_n, \mathcal
{S}(\bar{p}_n))$, and the above discussion shows that the
composite $\Psi_n \circ \hat{\psi}$ is the same homomorphism as that
which is functorially induced by the following commutative square
\begin{equation} \label{EAcomsq6}
\begin{CD}
\tilde {Y} \times_{\Gamma} \hat{\mathcal {K}}
@>>\subset>
\tilde{Y} \times_{\Gamma} \mathcal {E} \\
@VpVV @VVpV \\
\mathcal {K} @>>\subset> \mathcal {B}
\end{CD}
\end{equation}
whose horizontal arrows are both inclusions.  This fact together
with Remark \ref{Rfibration} shows immediately
that $\Psi$ is a surjection.
Hence to complete the proof that the top arrow in
diagram \eqref{EAcomlim} is an epimorphism,
it remains only to prove Proposition \ref{Pfinsubcomplex}.  For
this purpose we introduce some more notation.

\begin{Not}
If $X$ is a $\Gamma$-space and $T$ is a subgroup of
$\Gamma$, then
\[
X^T = \{x \in X \mid gx = x \text{ for all } g \in T\}.
\]
That is $X^T$ is the subspace fixed pointwise by $T$.
\end{Not}
\begin{Not}
For each pair of indices $n \ge m$, let
\[
\bar{f}_{mn} \colon \Gamma \times_{\Gamma_m} \mathcal {E}_m
\longrightarrow \Gamma \times_{\Gamma_n} \mathcal {E}_n
\]
be the $\Gamma$-equivariant map defined by the equation
$\bar{f}_{mn}([\gamma,x]) = [\gamma, f_{mn}(x)]$ for all $(\gamma,x)
\in \Gamma \times \mathcal {E}_m$.
\end{Not}

The key to proving Proposition \ref{Pfinsubcomplex}
is the following lemma.
\begin{Lem} \label{Lvircyclic}
For each virtually cyclic subgroup $T$ of $\Gamma$ and each
non-negative integer $i$
\[
\varinjlim_{n\in I} \pi_i \bigl((\Gamma \times_{\Gamma_n} \mathcal
{E}_n)^T\bigr) = 0.
\]
In particular $(\Gamma \times_{\Gamma_n} \mathcal {E}_n)^T$ is
nonempty for all sufficiently large indices $n$.
\end{Lem}

\begin{Rem} \label{Rdirsystem}
The maps in this directed system are induced by the family
$\bar{f}_{mn}$.  Note that, although base points need not be
preserved by the family $\bar{f}_{mn}$, this is irrelevant to the
vanishing assertion of the lemma.  In fact this assertion can be
paraphrased as saying that for any map
\[
\eta \colon S^i \longrightarrow (\Gamma
\times_{\Gamma_m} \mathcal {E}_m)^T,
\]
there exists an index $n \ge m$ such that
$\bar{f}_{mn} \circ \eta$ extends
to a map of $\mathbb {D}^{i+1}$ into $(\Gamma \times_{\Gamma_n}
\mathcal {E}_n)^T$.
\end{Rem}

The following fact is used in proving Lemma \ref{Lvircyclic}.
\begin{Rem} \label{Rvircyclic}
Let $S$ be a virtually cyclic subgroup of $\Gamma$.  Then there
exists an index $n \in I$ and a subgroup $\bar {S}$ of $\Gamma_n$
such that $\phi_n(\bar{S}) = S$ and
$\phi_n \colon \bar {S} \to S$ is
an isomorphism.  Such a group $\bar {S}$ is called a ``lift" of $S$
(to $\Gamma_n$).  Furthermore, if $\hat{S} \subseteq \Gamma_m$ is a
second lift, then there exists an index $s \ge m,n$ such that
$\phi_{ms}(\hat{S}) = \phi_{ns}(\bar {S})$.
\end{Rem}

The existence of lifts is a consequence of Remark \ref{Rfinpres}
because $S$ is finitely presented.  The ``uniqueness" is also easily
shown.

Here is another fact used in proving Lemma \ref{Lvircyclic}.

\begin{Rem} \label{Rcomsubset}
Given a compact subset $K$ of $\mathcal {E}_n$, there exists a
finitely generated subgroup $G$ of $\ker \phi_n$ such that
if $g \in \ker \phi_n$ and $gK \cap K \ne \emptyset$, then $g \in
G$.
\end{Rem}
\begin{proof}
We may assume that $K$ is a finite subcomplex of $\mathcal {E}_n$.
Let
\[
\mathcal {S} = \{ g\in \ker \phi_n \mid gK \cap K \ne \phi \}.
\]
For each pair of simplices $\sigma, \tau \in K$ such that $g\sigma =
\tau$ for some $g \in \mathcal {S}$, fix one choice of an element
$g_{\sigma \tau} \in \mathcal {S}$ such that $g_{\sigma \tau} \sigma
= \tau$.  Clearly the collection of such
choices $\{g_{\sigma \tau}\}$
forms a finite set.  Now given $g \in \mathcal {S}$, there exists a
pair of simplices $\sigma, \tau$ in $K$ such that $g \sigma = \tau$.
Therefore $g^{-1}g_{\sigma \tau}$ fixes $\sigma$ and is hence an
element in $\Gamma_n^{\sigma}$ which denotes the subgroup of
$\Gamma_n$ fixing $\sigma$.  But it is also in $\ker \phi_n$.
Therefore
\[
g^{-1}g_{\sigma \tau} \in \ker \phi_n \cap \Gamma_n^{\sigma}
\]
which is a finitely generated group contained in the set $\mathcal
{S}$.  Let $G$ be the group generated by $\{g_{\sigma \tau}\}$
together with the sets $\ker \phi_n \cap \Gamma_n^{\sigma}$ where
$\sigma$ varies over all simplices in $K$.  It is clearly a finitely
generated subgroup of $\ker \phi_n$ and contains $\mathcal {S}$.
\end{proof}

We list two more easily verified remarks needed to prove Lemma
\ref{Lvircyclic}.
\begin{Rem} \label{Requivclass}
If $[1,x] = [1,y]$ in $\Gamma
\times_{\Gamma_n} \mathcal {E}_n$, then $x =
\gamma y$ where $\gamma \in \ker \phi_n$.  (We use 1 to denote the
identity element of a group written multiplicatively.)  Consequently
if $[1, gx] = [1,x]$ for $g \in \Gamma_n$, then $g \in (\ker \phi_n)
\Gamma_n^x$.
\end{Rem}
\begin{Not}
The map $x \mapsto [1,x]$ is called the canonical projection of
$\mathcal {E}_n$ into $\Gamma \times_{\Gamma_n} \mathcal {E}_n$.
\end{Not}
\begin{Rem} \label{Rconsubset}
If $A$ is a connected subset of $\Gamma \times_{\Gamma_m} \mathcal
{E}_m$, then there exists an index $n \ge m$ such that $\bar{f}_{mn}
(A)$ is in the image of the canonical projection of $\mathcal {E}_n$
into $\Gamma \times_{\Gamma_n} \mathcal {E}_n$.
\end{Rem}

\begin{proof}[Proof of Lemma \ref{Lvircyclic}]
Let $\bar{T}$ be the lift of $T$ to $\Gamma_n$ given by Remark
\ref{Rvircyclic}.  Since $\mathcal{E}_n^{\bar{T}} \ne \emptyset$ and
its image under the canonical projection is contained in $(\Gamma
\times_{\Gamma_n} \mathcal {E}_n)^T$, we see that $(\Gamma
\times_{\Gamma_s} \mathcal {E}_s)^T \ne \emptyset$ for all
indices $s \ge n$.

Next let $x,y$ be two points in $(\Gamma
\times_{\Gamma_m} \mathcal
{E}_m)^T$.  Then by Remark \ref{Rconsubset} there are
points $\bar {x}, \bar{y} \in \mathcal {E}_n$ which canonically
project to $\bar{f}_{mn} (x) = \hat {x}$, $\bar {f}_{mn} (y) =
\hat{y}$ provided $n$ is sufficiently large.  By making $n$ perhaps
larger, we also have a lift $\bar {T}$ of $T$ to $\Gamma_n$.  If $g
\in \bar {T}$, then
\[
[1, g\bar{x}] = \phi_n(g) \hat{x} = \hat{x} = [1, \bar{x}]
\]
since $\hat{x} \in (\Gamma \times_{\Gamma_n}
\mathcal{E}_n)^T$.  Therefore by Remark \ref{Requivclass},
$g = \gamma \bar{g}$ where $\gamma \in \ker {\phi_n}$ and $\bar{g}
\bar {x} = \bar {x}$.  Consequently by picking a larger index $s$
such that $\phi_{ns} (\gamma) = 1$
\begin{equation} \label{Evircyclic1}
\tilde {g} \tilde {x} = \tilde{x}
\end{equation}
where $\tilde {g} = \phi_{ns} (g)$ and $\tilde {x} = f_{ns}
(\bar{x})$.  Note that $\tilde {g} \in \tilde{T} = \phi_{ns}(\bar
{T})$ which is a lift of $T$ to $\Gamma_s$.  Since $T$ is finitely
generated, we can in this way pick $s$ large enough so that both
equation \eqref{Evircyclic1} and
\begin{equation} \label{Evircyclic2}
\tilde {g} \tilde {y} = \tilde {y}
\end{equation}
hold for all $\tilde {g} \in \tilde {T}$,
where $\tilde{y}= f_{ns}(\bar
{y})$.  Hence $\tilde {x}, \tilde {y}$ are two points in the
contractible space $\mathcal {E}_s^{\tilde{T}}$ and consequently
connected by a path in $\mathcal {E}_s^{\tilde {T}}$.  The image of
this path under the canonical projection into $\Gamma
\times_{\Gamma_s} \mathcal{E}_s$ is the desired arc in $(\Gamma
\times_{\Gamma_s} \mathcal {E}_s)^T$ connecting $\bar{f}_{ms}(x)$ to
$\bar{f}_{ms}(y)$ and establishes the case of Lemma \ref{Lvircyclic}
when $i=0$.

Now consider the general case of a map
\[
\eta \colon S^i \longrightarrow (\Gamma \times_{\Gamma_m}
\mathcal {E}_m)^T
\]
cf.\ Remark \ref{Rdirsystem}.  We argue in a fashion similar to the
case $i=0$ but more complicated.  By Remark \ref{Rconsubset}, there
exists an index $n \ge m$ such that $\hat{\eta}
= \bar{f}_{mn} \circ \eta$
is in the image of the canonical projection
\[
\rho_n \colon \mathcal {E}_n \longrightarrow
\Gamma \times_{\Gamma_n} \mathcal
{E}_n.
\]
After a homotopy, we can assume that $\hat{\eta}$ is simplicial for
some simplicial structure on $S^i$ and is still in the image of
$\rho_n$ as well as inside $(\Gamma \times_{\Gamma_n} \mathcal
{E}_n)^T$.

For each simplex $\sigma$ of $S^i$, let $\hat{\eta}_{\sigma}$ denote
the restriction of $\hat {\eta}$ to $\sigma$ and let
$\bar{\eta}_{\sigma}$ be a lift of $\hat{\eta}_{\sigma}$ to
$\mathcal {E}_n$; i.e.\ $\rho_n \circ \bar{\eta}_{\sigma} =
\hat{\eta}_{\sigma}$ which exists since $\rho_n$ is a non-collapsing
simplicial map.  If $\tau$ is a face of $\sigma$, then
$\bar{\eta}_{\tau}$ may not be the restriction of
$\bar{\eta}_{\sigma}$ to $\tau$.
But by Remarks \ref{Rcomsubset} and
\ref{Requivclass}, there exists a larger index $s$ such that $\tilde
{\eta}_{\tau}$ is $\tilde {\eta}_{\sigma}$
restricted to $\tau$, where
\[
\tilde {\eta}_{\tau} = f_{ns}
\circ \bar{\eta}_{\tau} \quad \text{and}
\quad \tilde{\eta}_{\sigma} = f_{ns} \circ \bar{\eta}_{\sigma}.
\]
Hence the maps $\{\tilde {\eta}_{\sigma}\}$ where $\sigma$ varies
over the simplices in $S^i$ assemble to form a continuous map
$\tilde {\eta} \colon S^i \to \mathcal {E}_s$ such that
\[
\rho_s \circ \tilde {\eta} = \bar {f}_{ns} \circ \hat {\eta}.
\]
And Remark \ref{Rvircyclic} gives us a lift $\tilde{T}$ of $T$ to
$\Gamma_s$ provided $s$ was picked large enough.  Let $g \in \tilde
{T}$ and let $x \in S^i$.  Then
\[
[1, g\tilde{\eta}(x)] = \phi_s(g) \bar{f}_{ns}(\hat{\eta}(x)) =
\bar{f}_{ns} (\phi_s(g) \hat{\eta} (x)) =
\bar{f}_{ns} (\hat{\eta}(x)) = [1, \tilde{\eta} (x)]
\]
since $\phi_s(g) \in T$ and $\hat{\eta}(x) \in (\Gamma
\times_{\Gamma_n} \mathcal {E}_n)^T$.  Therefore by Remark
\ref{Requivclass}, $g = \gamma \tilde {g}$ where $\gamma \in \ker
\phi_s$ and $\tilde {g} \tilde {\eta} (x) = \tilde {\eta}(x)$.
Consequently by picking a larger index $t$ such
that $\phi_{st}(\gamma) =1$
\begin{equation} \label{Evircyclic3}
\dot{g} \dot{\eta} (x) = \dot{\eta}(x)
\end{equation}
where $\dot{g} = \phi_{st}(g)$ and $\dot{\eta} = f_{st} \circ
\tilde{\eta}$.  Note that $\dot{g} \in \dot {T} =
\phi_{st}(\tilde{T})$ which is a lift of $T$ to $\Gamma_t$.  Since
$T$ is finitely generated, Remark \ref{Rcomsubset} allows us to pick
$t$ large enough so that equation \eqref{Evircyclic3} holds for all
$\dot{g} \in \dot{T}$ and all $x \in S^i$.  Hence the image of
$\dot{\eta}$ lies in the contractible space $\mathcal
{E}_t^{\dot{T}}$ and consequently $\dot {\eta}$ extends to a map of
$\mathbb {D}^{i+1}$ into $\mathcal {E}_t^{\dot{T}}$.
The composition of this extended map with $\rho_t$ is an extension
of $\bar{f}_{nt} \circ
\hat{\eta}$ to a map of $\mathbb {D}^{i+1}$ into $(\Gamma
\times_{\Gamma_t} \mathcal {E}_t)^T$.  But
\[
\bar{f}_{mt} \circ \eta \sim \bar{f}_{nt} \circ \bar{f}_{mn} \circ
\eta \sim \bar{f}_{nt} \circ \hat{\eta}
\]
inside $(\Gamma \times_{\Gamma_t} \mathcal {E})^T$.  Therefore
$\bar{f}_{mt} \circ \eta$ also extends to a map of $\mathbb
{D}^{i+1}$ into $(\Gamma \times_{\Gamma_t} \mathcal {E}_t)^T$.
\end{proof}

\begin{proof}[Proof of Proposition \ref{Pfinsubcomplex}]
There are only a finite number of $\Gamma$-orbits of simplices in
$\hat{\mathcal {K}}$ since $\mathcal {K}$ is a finite complex.  Let
$\sigma_1, \sigma_2, \dots, \sigma_k$ be a complete irredundant list
of representatives for these $\Gamma$-orbits enumerated so that
$\dim \sigma_i \le \dim \sigma_j$ when $i \le j$.  Let $\hat
{\mathcal {K}}_i$ denote the subcomplex of $\hat{\mathcal {K}}$
defined by
\[
\hat{\mathcal {K}}_i = \bigcup_{j \le i} \Gamma \sigma_j.
\] Note that $\hat{\mathcal {K}_i}$ is a closed and
$\Gamma$-invariant subcomplex of $\hat{\mathcal {K}}$.  We will
construct $\psi$ by induction on $i$ over these subcomplexes; i.e.,
we assume that a $\Gamma$-equivariant map
\[
\psi_i \colon \hat{\mathcal {K}}_i \longrightarrow \Gamma
\times_{\Gamma_m} \mathcal {E}_m
\]
has already been constructed for some index $m \in I$ and then we
must construct a $\Gamma$-equivariant map
\[
\psi_{i+1} \colon \hat {\mathcal {K}}_{i+1} \longrightarrow \Gamma
\times_{\Gamma_n} \mathcal {E}_n
\]
for some index $n \ge m$.  Here is how this is done.  Note that
$\partial \sigma_{i+1}$ is contained in $\hat{\mathcal {K}}_i$ and
let
\[
\eta \colon \partial \sigma_{i+1} \longrightarrow \Gamma
\times_{\Gamma_m} \mathcal {E}_m
\]
be the restriction of $\psi_i$ to $\partial \sigma_{i+1}$.  Let $T =
\Gamma^{\sigma_{i+1}}$; i.e., $T$ is the subgroup of $\Gamma$ fixing
$\sigma_{i+1}$.  Then
\[
\im \eta \subseteq (\Gamma \times_{\Gamma_m} \mathcal {E}_m)^T
\]
since $\eta$ is $\Gamma$-equivariant and $T$ fixes $\partial
\sigma_{i+1}$.  Since $\partial \sigma_{i+1}$ is homeomorphic to a
sphere and $\sigma_{i+1}$ to the ball that sphere bounds, Lemma
\ref{Lvircyclic} yields that there exists an index $n \ge m$ such
that $\bar{f}_{mn} \circ \eta$ extends to a map
\[
\bar{\eta} \colon \sigma_{i+1} \longrightarrow
(\Gamma \times_{\Gamma_n} \mathcal {E}_n)^T
\]
(see Remark \ref{Rdirsystem}).  Now we define $\psi_{i+1}$ by
\[
\psi_{i+1}(x)=
\begin{cases}
\bar{f}_{mn}(\psi_i(x)) &\text{if } x \in \hat{\mathcal {K}_i} \\
\gamma \bar{\eta}(\gamma^{-1} x) &\text{if } x \in \gamma
\sigma_{i+1}.
\end{cases}
\]
This map is easily seen to be well-defined, continuous and
$\Gamma$-equivariant and the proof of Proposition
\ref{Pfinsubcomplex} is complete.
\end{proof}

\begin{Rem} \label{Raddendum}
A version of the FIC can be formulated for any full class
$\mathcal {C}$ of subgroups of $\Gamma$; i.e., one closed with
respect to subgroups and conjugates (see \cite[introduction and
appendix]{farrelljonesiso}).  We say
that the full class $\mathcal {C}$ is \emph{small}
if it satisfies the following two additional properties:
\begin{enumerate}
\item Each group in $\mathcal {C}$ is finitely presented.
\item The class $\mathcal {C}$ is closed with respect to quotient
groups.
\end{enumerate}
Some examples of small classes are:
\begin{enumerate}
\item all virtually cyclic groups
\item all finite groups
\item the trivial group.
\end{enumerate}
\end{Rem}

The proof given above for Proposition \ref{Pfinsubcomplex} actually
remains valid when the class of virtually cyclic groups is replaced
by any small class of groups.  In particular, when it is replaced by
the class consisting of only the trivial group, the above proof
yields Addendum \ref{A}, which we state after defining the following
notation.

\begin{Not}
Given a subspace $K$ of $B$, let $\hat {K}$ denote its inverse image
under the projection $E \to B$.
\end{Not}

\begin{Add} \label{A}
Given A finite subcomplex $K$ of $B$, there exists an index $n \in I$
and a $\Gamma$-equivariant map $\psi \colon \hat{K} \to \Gamma
\times_{\Gamma_n} E_n$.
\end{Add}

We now develop some more ideas needed to show that $\pi_i(\mathfrak
{a})$ is monic.

\begin{Not}
Given an index $n \in  I$ and a subspace $K_n \subseteq B_n$, let
$\hat{K}_n$ denote its inverse image under the projection $E_n \to
B_n$.  Similarly if $K\subseteq B$ or $K_n \subseteq B_n$, define
$\check{K}$ and $\check{K}_n$ to be their inverse images under the
projections $Y \to B$ and $Y_n \to B_n$, respectively.  Denote by
$\tilde {K}$ and $\tilde{K}_n$ the spaces in the upper left corners
of the Cartesian squares
\begin{equation} \label{EAcomsq11}
\begin{CD}
\tilde{K} @>>> \hat{K} \\
@VVV @VVV \\
\check{K} @>>> K
\end{CD}
\end{equation}

\begin{equation} \label{EAcomsq12}
\begin{CD}
\tilde{K}_n @>>> \hat{K}_n \\
@VVV @VVV \\
\check{K}_n @>>> K_n.
\end{CD}
\end{equation}

Notice that $\tilde{K} \to \check{K}$ and $\tilde {K}_n \to
\check{K}_n$ are regular $\Gamma$ and $\Gamma_n$-covering spaces,
respectively.

Now fix a finite subcomplex $K$ of $B$.  Because of Addendum
\ref{A}, there exists an index $n \in I$ and a $\Gamma$-equivariant
map $\psi \colon \hat{K} \to \Gamma
\times_{\Gamma_n} E_n$.  Let $\bar{\psi}
\colon K \to B_n$ denote the map induced on orbit spaces by $\psi$.
Also define a $\Gamma$-equivariant map
\[
\bar{f}_n \colon \Gamma \times_{\Gamma_n} E_n \longrightarrow E
\]
as in the paragraph following Proposition \ref{Pfinsubcomplex}.
Then $\bar{f}_n$
induces $\bar{\phi}_n \colon B_n \to B$ on orbit spaces (note that
when the small class consists of only the trivial group, then
$\mathcal {E} = E$ and $\mathcal {E}_n = E_n$) and the composite
$\bar{f}_n \circ \psi \colon \hat{K} \to E$ is $\Gamma$-equivariantly
homotopic to the inclusion $\hat{K} \hookrightarrow E$.  Passing to
orbit spaces, this homotopy determines a homotopy $h_t$ of
$\bar{\phi}_n \circ \bar{\psi}$ to the inclusion $K \hookrightarrow
B$.
\end{Not}

\begin{Lem} \label{L6.16}
There exists a bundle map $\check{\psi} \colon \check{K} \to Y_n$
covering $\bar{\psi}$ and a homotopy of $\hat{\phi}_n \circ \check
{\psi}$ to the inclusion $\check{K} \hookrightarrow Y$ via bundle
maps covering $h_t$.
\end{Lem}
\begin{proof}
This is a consequence of the covering homotopy theorem.
\end{proof}

By considering the Cartesian square
\begin{equation} \label{EAcomsq13}
\begin{CD}
\Gamma \times _{\Gamma_n} \tilde{Y}_n @>>>
\Gamma \times_{\Gamma_n} E_n \\
@VVV @VVV \\
Y_n @>>> B_n
\end{CD}
\end{equation}
together with square \eqref{EAcomsq11}, it is seen that the triple
of maps $\bar{\psi} \colon K \to B_n$, $\psi \colon \hat{K} \to
\Gamma \times_{\Gamma_n} E_n$, $\check{\psi} \colon \check{K} \to
Y_n$ determine a $\Gamma$-equivariant map
\[
\tilde{\psi} \colon \tilde{K} \longrightarrow \Gamma
\times_{\Gamma_n} \tilde{Y}_n.
\]
Also by considering the Cartesian square
\begin{equation} \label{EAcomsq14}
\begin{CD}
\tilde{Y} @>>> E \\
@VVV @VVV \\
Y @>>> B
\end{CD}
\end{equation}
together with square \eqref{EAcomsq13}, it is seen that the triple
$\bar{\phi}_n \colon B_n \to B$, $\bar{f}_n \colon \Gamma
\times_{\Gamma_n} E_n \to E$, $\hat{\phi}_n \colon Y_n \to Y$
determines a $\Gamma$-equivariant map
\[
\tilde{f}_n \colon \Gamma \times_{\Gamma_n} \tilde{Y}_n
\to \tilde{Y}.
\]
(In fact $\tilde{f}_n$ is determined from $\tilde{\phi}_n \colon
\tilde{Y}_n \to \tilde{Y}$ by the equation $\tilde{f}_n([\gamma,x]) =
\gamma \tilde{\phi}_n(x)$ for all $(\gamma,x) \in \Gamma \times
\tilde{Y}_n$.)
\begin{Lem} \label{L6.17}
The map $\tilde{f}_n \circ \tilde{\psi} \colon \tilde{K} \to
\tilde{Y}$ is $\Gamma$-equivariantly homotopic to the inclusion map
$\tilde {K} \hookrightarrow \tilde{Y}$.
\end{Lem}
\begin{proof}
This follows from Lemma \ref{L6.16} since $\tilde{f}_n \circ
\tilde{\psi}$ is determined from the Cartesian squares
\eqref{EAcomsq11} and \eqref{EAcomsq14} by the triple of maps
$\bar{\phi}_n \circ \bar{\psi} \colon K \to B$,
$\bar{f}_n \circ \psi \colon \hat{K} \to E$ and $\hat{\phi}_n \circ
\check{\psi} \to Y$.
\end{proof}

The $\Gamma$-spaces $\tilde{K}$ and $\Gamma \times_{\Gamma_n}
\tilde{Y}_n$ determine two new continuous maps onto
$\mathcal {B}$, namely
\[
p_K \colon \tilde{K} \times _{\Gamma} \mathcal{E} \to
\mathcal{E}/\Gamma = \mathcal {B} \quad\text{and}\quad
p_n' \colon (\Gamma \times_{\Gamma_n} \tilde{Y}_n) \times_{\Gamma}
\mathcal {E} \to \mathcal {E}/\Gamma = \mathcal{B}.
\]
Furthermore, the $\Gamma$-equivariant maps
$\tilde {\psi} \colon \tilde{K} \to \Gamma \times_{\Gamma_n}
\tilde{Y}_n$, $\tilde{f}_n \colon
\Gamma \times_{\Gamma_n} \tilde{Y}_n \to
\tilde{Y}$, and the inclusion map $\tilde{K} \hookrightarrow
\tilde{Y}$ determine group homomorphisms
\begin{align*}
\alpha &\colon \mathbb {H}_i(\mathcal{B}, \mathcal{S}(p_K))
\longrightarrow
\mathbb{H}_i(\mathcal {B}, \mathcal{S}(p_n')), \\
\beta &\colon \mathbb {H}_i(\mathcal {B}, \mathcal {S}(p_n'))
\longrightarrow
\mathbb {H}_i(\mathcal {B}, \mathcal {S}(p)), \\
\gamma &\colon \mathbb {H}_i (\mathcal {B}, \mathcal {S}(p_K))
\longrightarrow
\mathbb {H}_i(\mathcal {B}, \mathcal {S}(p)),
\end{align*}
respectively.

\begin{Lem} \label{L6.18}
$\gamma = \beta \circ \alpha$.
\end{Lem}
\begin{proof}
This follows directly from Lemma \ref{L6.17}.
\end{proof}

Now fix an index $m \in I$ and a finite subcomplex $K_m \subseteq
B_m$ such that $\bar{\phi}_m (K_m) \subseteq K$.  Note that
$\hat{\phi}_m$ restricted to $\check{K}_m$ is a bundle map into
$\check{K}$ which covers the restriction of $\bar {\phi}_m$ to $K_m$.
Also the functorially defined $\phi_m$-equivariant map $E_m \to E$
restricts to a map $\hat{K}_m \to \hat{K}$ covering
$\bar{\phi}_m$.  Using the Cartesian squares \eqref{EAcomsq12}, in
which $m$ replaces $n$, and \eqref{EAcomsq11}, we see that the triple
of maps $\bar{\phi}_m \colon K_m \to K$, $\hat{K}_m \to K$, and
$\hat{\phi}_m \colon \check{K}_m \to K$ determine a
$\phi_m$-equivariant map $\tilde{K}_m \to \tilde {K}$ which in turn
determines a group homomorphism
\[
\delta \colon \mathbb {H}_i(\mathcal {B}_m, \mathcal {S}(p_{K_m}))
\longrightarrow \mathbb {H}_i(\mathcal {B}, \mathcal {S}(p_K))
\]
where the surjection 
\[
p_{K_m} \colon \tilde {K}_m \times_{\Gamma_m} \mathcal {E}_m
\longrightarrow \mathcal {E}_m/\Gamma_m = \mathcal {B}_m
\]
comes from projection onto the second factor of $\tilde{K}_m \times
\mathcal {E}_m$.  Also the Cartesian product of the inclusion $\tilde
{K}_m \hookrightarrow \tilde{Y}_m$ with $\id_{\mathcal {E}_m}$
determines a group homomorphism
\[
\epsilon \colon \mathbb {H}_i(\mathcal {B}_m, \mathcal {S}(p_{K_m}))
\longrightarrow \mathbb {H}_i (\mathcal {B}_m, \mathcal {S}(p_m)).
\]
Let
\[
\omega_m \colon \mathbb {H}_i(\mathcal {B}_m, \mathcal {S}(p_m))
\longrightarrow \mathbb {H}_i(\mathcal {B}, \mathcal {S}(p))
\]
denote the homomorphism $\Psi _m \circ \Phi_m$ (see the paragraph
containing diagram \eqref{EAcomsq1}).

\begin{Lem} \label{L6.19}
$\omega_m \circ \epsilon = \gamma \circ \delta$.
\end{Lem}
\begin{proof}
Consider the $\phi_m$-equivariant map determined by the Cartesian
square \eqref{EAcomsq12}, in which $m$ replaces $n$, and
\eqref{EAcomsq14} via the triple of continuous maps $\bar{\phi}_m$,
the functorial map $E_m \to E$, and $\hat{\phi}_m$, where these three
maps are restricted to $K_m$, $\hat{K}_m$ and $\check{K}_m$,
respectively.  Then both homomorphisms $\omega_m \circ \epsilon$ and
$\gamma \circ \delta$ are induced by the Cartesian product of this
map with the $\phi_m$-equivariant map $f_m \colon \mathcal {E}_m \to
\mathcal {E}$.
\end{proof}

Now notice that the universal $(\Gamma, \mathcal {C})$-space
$\mathcal {E}$ is also a $\Gamma_n$-space via the homomorphism
$\phi_n \colon \Gamma_n \to \Gamma$.  Define a new full class of
subgroups $\mathcal {C}_n$ of $\Gamma_n$ by
\[
\mathcal {C}_n = \{\text{subgroups $S$ of } \Gamma_n \mid \phi_n(S)
\in \mathcal {C} \}.
\]
Then the following fact is an immediate consequence of the
definitions.

\begin{Rem} \label{R6.20}
The $\Gamma_n$-space $\mathcal {E}$ is universal for the class
$\mathcal {C}_n$.
\end{Rem}

Denote by $\mathcal {B}_n^*$ the orbit space $\mathcal {E}/\Gamma_n$.
Then projection onto the second factor of $\tilde {Y}_n \times
\mathcal {E}$ defines as usual a continuous surjection 
\[
p_n^* \colon \tilde{Y}_n \times _{\Gamma_n} \mathcal {E}
\longrightarrow \mathcal {E}/\Gamma_n = \mathcal {B}_n^*.
\]
Note that there is a natural identification
\[
\tilde{Y}_n \times_{\Gamma_n} \mathcal {E} = (\Gamma
\times_{\Gamma_n} \tilde{Y}_n) \times_{\Gamma} \mathcal {E}
\]
and that $p_n'$ factors through $p_n^*$; i.e., we have the following
commutative diagram:
\begin{equation}  \label{EAcomsq15}
\begin{CD}
\tilde{Y}_n \times _{\Gamma_n} \mathcal {E} @= (\Gamma
\times_{\Gamma_n} \tilde {Y}_n) \times_{\Gamma } \mathcal {E} \\
@Vp_n^*VV @VVp_n'V \\
\mathcal{B}_n^* = \mathcal {E}/\Gamma_n @>>> \mathcal {E}/\Gamma =
\mathcal {B}
\end{CD}
\end{equation}
where $\mathcal {E}/\Gamma_n \to \mathcal {E}/\Gamma$
is the canonical quotient map.  Let
\[
\Theta \colon \mathbb {H}_i(\mathcal {B}_n^*, \mathcal {S}(p_n^*))
\longrightarrow \mathbb {H}_i(\mathcal {B}, \mathcal {S}(p_n'))
\]
be the homomorphism induced by commutative square \eqref{EAcomsq15}.
If $\phi_n \colon \Gamma_n \to \Gamma$ is a surjection, then the
canonical quotient map $\mathcal{E}/\Gamma_n \to \mathcal
{E}/\Gamma$ is the identity and hence so is $\Theta$.  In the general
case, $\mathcal {E}/\Gamma_n \to \mathcal {E}/\Gamma$ is a simplicial
map which collapses no simplices.  The following remark is a
consequence of this observation and \cite[Addendum A.6.1 and Lemma
A.7]{farrelljonesiso}.
\begin{Rem} \label{R6.21}
$\Theta$ is an isomorphism.
\end{Rem}
The next remark follows directly from \cite[Lemma A.12.1 and Theorem
A.8]{farrelljonesiso}.

\begin{Rem} \label{R6.22}
Let $C$ and $C'$ be two full classes of subgroups of a group $\pi$
with $C' \subseteq C$.  If $\pi$ satisfies the FIC with respect to
$C'$ (cf.\ Remark \ref{Raddendum}), then $\pi$ also satisfies the
FIC with respect to $C$.
\end{Rem}

\begin{Lem} \label{L6.23}
The assembly map $\mathbb {H}_i( \mathcal {B}_n^*, \mathcal
{S}(p_n^*)) \to \pi_i(\mathcal {S}(Y_n))$ is an isomorphism.
\end{Lem}
\begin{proof}
Since $\Gamma_n$ satisfies the FIC with respect to $\mathcal {C}$ by
assumption and $\mathcal {C}$ is a subclass of $\mathcal{C}_n$, we
see that $\Gamma_n$ also satisfies the FIC with respect to
$\mathcal {C}_n$.  This is because of Remark \ref{R6.22} in which we
set $\pi = \Gamma_n$, $C' = \mathcal {C}$ and $C = \mathcal {C}_n$.
The Lemma follows immediately from this fact.
\end{proof}

We are now ready to complete the proof of Theorem \ref{Tdirlimit} by
showing that the map $\pi_1(\mathfrak{a})$ in diagram \eqref{EAcomlim}
is monic.  For this purpose, let $a \in \mathbb {H}_i(\mathcal {B},
\mathcal {S}(p))$ which ``assembles" to 0 in $\pi_i(\mathcal
{S}(Y))$; i.e., such that
\[
\pi_i(\mathfrak{a})(a) = 0.
\]
We need to show that $a= 0$.  Since the top arrow in diagram
\eqref{EAcomlim} is an epimorphism, there exists an index $m\in I$ and
an element $a_1 \in \mathbb {H}_i(\mathcal {B}_m, \mathcal
{S}(p_m))$ which maps to $a$ via $\omega_n$; i.e.,
\begin{equation} \label{E6.16}
\omega_m(a_1) = a.
\end{equation}
Then $a_1$ assembles to
an element $b_1 \in \pi_i(\mathcal {S}(Y_m))$.
Because of Remark \ref{Rfibration}, there exists a finite subcomplex
$K_m$ of $B_m$ and an element $b_2 \in \pi_i(\mathcal
{S}(\check{K}_m))$ which maps to $b_1$ under the map induced by the
inclusion $\check{K}_m \hookrightarrow Y_m$.  Since $b_1$ maps to 0
in $\pi_i(\mathcal {S}(Y))$ under the map induced by $\hat{\phi}_m$,
Remark \ref{Rfibration} also shows that there exists a finite
subcomplex $K$ of $B$ such that both $\bar{\phi}_m(K_m) \subseteq K$
and $b_2$ maps to 0 in
$\pi_i(\mathcal {S}(\check{K}))$ under the map
induced by the restriction of $\hat{\phi}_m$ to $\check{K}_m$.  Now
$b_2$ is in the image of the assembly map because $\Gamma_m$
satisfies the FIC; i.e., there is an element $a_2 \in \mathbb
{H}_i(\mathcal {B}_m, \mathcal {S}(p_{K_m}))$ which assembles to
$b_2$.  Note that
\begin{equation} \label{E6.17}
\epsilon(a_2) = a_1
\end{equation}
since $\epsilon (a_2)$ and $a_1$ both assemble to $b_1$ and
$\Gamma_m$ satisfies the FIC.  Therefore
\begin{equation} \label{E6.18}
\gamma(\delta(a_2)) = a
\end{equation}
because of equations \eqref{E6.16}, \eqref{E6.17} and Lemma
\ref{L6.19}.  Also $\delta(a_2)$ assembles to 0 in $\pi_i(\mathcal
{S}(\check{K}))$ since $b_2$ maps to 0.  Consequently
\begin{equation} \label{E6.19}
a_3 = \alpha(\delta(a_2))
\end{equation}
assembles to 0 in $\pi_i(\mathcal {S}(Y_n))$.  Note there is a
natural identification
\[
Y_n = (\Gamma \times_{\Gamma_n} \tilde{Y}_n)/\Gamma.
\]
Furthermore, notice that
\begin{equation} \label{E6.20}
\beta(a_3) = a
\end{equation}
because of equations \eqref{E6.18}, \eqref{E6.19} and Lemma
\ref{L6.18}.  Using Remark \ref{R6.21}, we find there is an element
$a_4 \in \mathbb {H}_i(\mathcal {B}_n^*, \mathcal {S}(p_n^*))$ such
that
\begin{equation} \label{E6.21}
\Theta(a_4) = a_3.
\end{equation}
Furthermore, $a_4$ assembles to 0 in $\pi_i(\mathcal {S}(Y_n))$
since $a_3$ does.  Therefore
\begin{equation} \label{E6.22}
a_4 = 0
\end{equation}
because of Lemma \ref{L6.23}.  Concatenating equations \eqref{E6.20},
\eqref{E6.21}, \eqref{E6.22} shows that $a=0$.  This completes the
proof of Theorem \ref{Tdirlimit}.
\end{proof}

The proof given above of Theorem \ref{Tdirlimit} actually proved the
following more general statement.
\begin{Add}
Let $\Gamma_n$, $n \in I$, be a directed system of groups with
\[
\Gamma = \varinjlim_{n\in I} \Gamma_n
\]
and let $\mathcal {C}$ be a small class of groups.  If each group
$\Gamma_n$ satisfies the FIC with respect to $\mathcal {C}$,
then so does $\Gamma$.
\end{Add}

\bibliographystyle{plain}

\begin{thebibliography}{10}

\bibitem{bass}
Hyman Bass, \emph{Algebraic ${K}$-theory}, W. A. Benjamin, Inc., New
  York-Amsterdam, 1968.

\bibitem{baumslagdyer}
Gilbert Baumslag and Eldon Dyer, \emph{The integral homology of finitely
  generated metabelian groups. {I}}, Amer. J. Math. \textbf{104} (1982), no.~1,
  173--182.

\bibitem{brown}
Kenneth~S. Brown, \emph{Cohomology of groups}, Springer-Verlag, New York, 1994,
  Corrected reprint of the 1982 original.

\bibitem{BrownGeoghegan84}
Kenneth~S. Brown and Ross Geoghegan, \emph{An infinite-dimensional torsion-free
  {$FP_{\infty}$} group}, Invent. Math. \textbf{77} (1984), no.~2, 367--381.

\bibitem{Cappell74}
Sylvain~E. Cappell, \emph{On connected sums of manifolds}, Topology \textbf{13}
  (1974), 395--400.

\bibitem{dicksdunwoody}
Warren Dicks and M.~J. Dunwoody, \emph{Groups acting on graphs}, Cambridge
  University Press, Cambridge, 1989.

\bibitem{dixonmortimer}
John~D. Dixon and Brian Mortimer, \emph{Permutation groups}, Springer-Verlag,
  New York, 1996.

\bibitem{Farley00}
Dan Farley, Ph.D. thesis, SUNY Binghamton, 2000.

\bibitem{farrellhsiangformula}
F.~T. Farrell and W.-C. Hsiang, \emph{A formula for ${K}\sb{1}{R}\sb{\alpha
  }\,[{T}]$}, Applications of Categorical Algebra (Proc. Sympos. Pure Math.,
  Vol. XVII, New York, 1968), Amer. Math. Soc., Providence, R.I., 1970,
  pp.~192--218.

\bibitem{farrellhsiangjlms}
F.~T. Farrell and W.~C. Hsiang, \emph{The {W}hitehead group of poly-(finite or
  cyclic) groups}, J. London Math. Soc. (2) \textbf{24} (1981), no.~2,
  308--324.

\bibitem{farrelljonesiso}
F.~T. Farrell and L.~E. Jones, \emph{Isomorphism conjectures in algebraic
  ${K}$-theory}, J. Amer. Math. Soc. \textbf{6} (1993), no.~2, 249--297.

\bibitem{farrellroushon}
F.~T. Farrell and Sayed~K. Roushon, \emph{The {W}hitehead groups of braid
  groups vanish}, Internat. Math. Res. Notices (2000), no.~10, 515--526.

\bibitem{klm}
P.~H. Kropholler, P.~A. Linnell, and J.~A. Moody, \emph{Applications of a new
  ${K}$-theoretic theorem to soluble group rings}, Proc. Amer. Math. Soc.
  \textbf{104} (1988), no.~3, 675--684.

\bibitem{quinninv}
Frank Quinn, \emph{Ends of maps. {I}{I}}, Invent. Math. \textbf{68} (1982),
  no.~3, 353--424.

\bibitem{ranicki}
Andrew Ranicki, \emph{Exact sequences in the algebraic theory of surgery},
  Princeton University Press, Princeton, N.J., 1981.

\bibitem{robinson}
Derek J.~S. Robinson, \emph{A course in the theory of groups}, second ed.,
  Springer-Verlag, New York, 1996.

\bibitem{scottwall}
Peter Scott and Terry Wall, \emph{Topological methods in group theory},
  Homological group theory (Proc. Sympos., Durham, 1977), Cambridge Univ.
  Press, Cambridge, 1979, pp.~137--203.

\bibitem{waldhausen}
Friedhelm Waldhausen, \emph{Algebraic ${K}$-theory of generalized free
  products. {I}{I}{I}, {I}{V}}, Ann. of Math. (2) \textbf{108} (1978), no.~2,
  205--256.

\bibitem{Wall76}
C.~T.~C. Wall, \emph{Classification of {H}ermitian {F}orms. {V}{I}. {G}roup
  rings}, Ann. of Math. (2) \textbf{103} (1976), no.~1, 1--80.

\bibitem{weinberger}
Shmuel Weinberger, \emph{There exist finitely presented groups with infinite
  ooze}, Duke Math. J. \textbf{49} (1982), no.~4, 1129--1133.

\end{thebibliography}

\providecommand{\bysame}{\leavevmode\hbox to3em{\hrulefill}\thinspace}

\end{document}